\documentclass[a4paper,twoside,10pt]{amsart}
\usepackage{amsmath,amsfonts,amssymb,amsthm}
\newtheorem{theorem}{Theorem}[section]
\newtheorem{proposition}[theorem]{Proposition}
\usepackage{graphicx}
\usepackage[colorlinks=true,pdfstartview=FitH,pdfview=FitH]{hyperref}
\numberwithin{equation}{section}

\def\placefigure#1#2{%
\setbox0=\hbox{\includegraphics[width=#2]{#1}}%
\hbox{\copy0}\hglue -\wd0}%

\def\placetext#1#2#3{%
\setbox1=\hbox{#1}\hglue#2\wd0\lower-#3\wd0\copy1\hglue-\wd1\hglue-#2\wd0}%

\def\im{\mathrm{i}}
\def\d{\mathrm{d}}
\def\e{\mathrm{e}}
\def\erfc{\operatorname{erfc}}

\begin{document}

\title[Asymptotics for the incomplete gamma function]{Asymptotic expansions for the incomplete gamma function in the transition regions}

\author[G. Nemes]{Gerg\H{o} Nemes}
\email{Gergo.Nemes@ed.ac.uk}
\author[A. B. Olde Daalhuis]{Adri B. Olde Daalhuis}
\email{A.B.Olde.Daalhuis@ed.ac.uk}
\address{School of Mathematics, The University of Edinburgh, James Clerk Maxwell Building, The King's Buildings, Peter Guthrie Tait Road, Edinburgh EH9 3FD, UK}
\thanks{The authors' research was supported by a research grant (GRANT11863412/70NANB15H221) from the National Institute of Standards and Technology.}


\begin{abstract}
We construct asymptotic expansions for the normalised incomplete gamma function $Q(a,z)=\Gamma(a,z)/\Gamma(a)$ that are valid in the transition regions, including the case $z\approx a$,  and have simple polynomial coefficients. For Bessel functions, these type of expansions are well known, but for the normalised incomplete gamma function they were missing from the literature. A detailed historical overview is included. We also derive an asymptotic expansion for the corresponding inverse problem, which has importance in probability theory and mathematical statistics. The coefficients in this expansion are again simple polynomials, and therefore its implementation is straightforward.
As a byproduct, we give the first complete asymptotic expansion as $a\to-\infty$ of the unique negative zero of the regularised incomplete gamma function $\gamma^*(a,x)$.
\end{abstract}

\maketitle


\section{Introduction and main results}\label{section0}
The normalised incomplete gamma function $Q(a,z)=\Gamma(a,z)/\Gamma(a)$ is one of the most widely used special functions of two variables. It is used in constructing gamma distributions \cite[Ch. 17]{JKB94}, which appear naturally in the theory associated with normally distributed random variables. In many applications, Poisson random variables are used in which the Poisson rate is not fixed, and the function $Q(a,z)$ plays an important role in these cases (see, for example, \cite{Giles2016}). Consequently, efficient and accurate approximations for the incomplete gamma function are essential. In many applications the inverse problem also arises naturally, that is, the problem of determining the solution $x$, real and positive, of the equation $Q(a,x)=q$ for given $0<q<1$ and $a>0$.

In the important papers \cite{Temme1979, Temme1996b}, Temme established asymptotic expansions for the normalised incomplete gamma function as $a\to\infty$, that are valid uniformly for $|z|\in[0,\infty)$. The asymptotic expansions of the corresponding inverse functions are discussed in \cite{Temme1992}. The large region of validity is important in applications, and many papers in physics, applied statistics, network and control theory, and so on, refer to these uniform asymptotic expansions (see, e.g., \cite{Lukas2014, Navarra2010, Paillard2008, Palffy2017, Ravelomanana2004, Trailovic2004}).

However, applying Temme's expansions in the neighbourhood of the (most interesting) point $z = a$ can be difficult because its coefficients possess a removable singularity at this point. To overcome this difficulty, Temme \cite{Temme1979} gave power series expansions for these coefficients. Having to use local power series expansions for the coefficients in uniform asymptotic expansions, is not a very elegant solution. In this paper, we propose an alternative approach to the problem. 

We will construct what we call \emph{transition region expansions}, and provide full details of the inversion of these new expansions. 
These are expansions that are valid in the regions in which $Q(a,z)$ changes dramatically, and their coefficients are polynomials satisfying simple recurrence relations. The region of validity overlaps with those of the non-uniform ``outer'' expansions. Furthermore, the coefficients of their inversions are simple polynomials, whose computation and implementation are straightforward.

It is surprising to us that these transition region expansions for the normalised incomplete gamma function have not yet been discussed in the literature, given the fact that expansions of similar type for Bessel functions are well known (see \cite[\href{http://dlmf.nist.gov/10.19.iii}{\S 10.19(iii)}]{NIST:DLMF}). What our new expansions and the transition region expansions for the Bessel functions have in common is that both mimic the corresponding uniform expansions; compare the similarities between \eqref{eq7} and \eqref{eq5}, and between \cite[\href{http://dlmf.nist.gov/10.19.E8}{Eq. 10.19.8}]{NIST:DLMF} and \cite[\href{http://dlmf.nist.gov/10.20.E4}{Eq. 10.20.4}]{NIST:DLMF}. The ideas in this paper should also be applicable to other cumulative distribution functions using the results of \cite{Temme1982}.

The proofs of our new results rely heavily on the uniform asymptotic expansions by Temme. Thus, we shall provide full details of his expansions, and even include new recurrence relations for the coefficients in the local power series expansions mentioned above (see Appendix \ref{appendixb}). For more information about existing results in the literature and their connection with those we prove in this paper, the reader is referred to the historical overview in Section \ref{section1}.

The following theorem is the main result of the paper.

\begin{theorem}\label{thm1} The normalised incomplete gamma function admits the asymptotic expansions
\begin{equation}\label{eq7}
Q\left(a,a+\tau a^{\frac{1}{2}}\right) \sim \tfrac{1}{2} \erfc\left(2^{-\frac{1}{2}}\tau\right) + \frac{1}{\sqrt {2\pi a}} 
\exp\left(-\frac{\tau^2}{2}\right) \sum_{n=0}^\infty\frac{C_n(\tau)}{a^{n/2}}
\end{equation}
and
\begin{gather}\label{eq8}
\begin{split}
\frac{\e^{\pm\pi \im  a}}{2\im  \sin (\pi a)}Q\left(-a,(a+\tau a^{\frac{1}{2}})\e^{\pm\pi \im  }\right) \sim & \pm \tfrac{1}{2}\erfc 
\left( \pm 2^{-\frac{1}{2}} \im  \tau \right) 
\\  &- \frac{\im   }{\sqrt{2\pi a}}\exp \left(\frac{\tau^2}{2}\right)\sum_{n=0}^\infty\frac{(-\im   )^n C_n(\im   \tau)}{a^{n/2}},
\end{split}
\end{gather}
as $a\to \infty$ in the sector $|\arg a| \le \pi-\delta<\pi$, uniformly with respect to bounded complex values of $\tau$. Here, $\erfc$ denotes the complementary error function \cite[\href{http://dlmf.nist.gov/7.2.E2}{Eq. 7.2.2}]{NIST:DLMF}. The coefficients $C_n(\tau)$ are polynomials in $\tau$ of degree $3n+2$ and satisfy
\[
C_0 (\tau ) = \frac{1}{3}\tau^2-\frac{1}{3},
\]
\[
C_n (\tau ) + \tau C'_n (\tau ) - C''_n (\tau ) = \tau (\tau ^2  - 2)C_{n - 1} (\tau ) - (2\tau^2  - 1)C'_{n - 1} (\tau ) + \tau C''_{n - 1} (\tau )
\]
for $n\geq 1$. In addition, the even- and odd-order polynomials are even and odd functions, respectively.
\end{theorem}

The explicit forms of the next two coefficients in \eqref{eq7} are as follows:
\[
C_1 (\tau) = \frac{1}{18}\tau^5-\frac{11}{36}\tau^3+\frac{1}{12}\tau,\quad C_2 (\tau) = \frac{1}{162}\tau^8-\frac{29}{324}\tau^6+\frac{133}{540}\tau^4-\frac{23}{540}\tau^2-\frac{1}{540}.
\]
For higher coefficients and for an algorithm to generate them, see the end of Section \ref{section2}.

The most common case in applications is when both variables of $Q(a,z)$ are real and positive. 
It is possible to combine our transition region expansion \eqref{eq7}, with the outer expansions that exist in the literature and cover the 
complete range $\tau > -a^{1/2}$  (or equivalently, the range $z>0$). First we note that in Theorem \ref{thm1} we assume that $\tau$ is bounded. 
We can weaken this restriction as follows. 
Since $C_n(\tau)$ is a polynomial of degree $3n+2$ it follows that we can allow $\tau = \mathcal{O}\left( |a|^\mu\right)$ as $a\to \infty$, 
as long as $\mu < \frac{1}{6}$. See the proposition below. In Section \ref{section1}, we will give the expansions \eqref{eq1} and \eqref{eq2}. These are also asymptotic expansions with simple coefficients, and with the notation of this section, these two outer expansions are valid for $|\tau| \geq a^\mu$ as $a\to +\infty$, with any $\mu >0$ (see Appendix \ref{appendixa}). Therefore, by combining these three expansions, we can cover the whole range $\tau > -a^{1/2}$.

\begin{proposition}\label{prop1}
The expansions \eqref{eq7} and \eqref{eq8} are still valid in the sense of generalised asymptotic expansions (see, e.g., 
\cite[\href{http://dlmf.nist.gov/2.1.v}{\S 2.1(v)}]{NIST:DLMF})
provided that $\tau = \mathcal{O}\left( |a|^\mu\right)$ as $a\to \infty$ with some constant $\mu < \frac{1}{6}$.
\end{proposition}

In the inverse problem we are trying to find $x=x(a,q)$ such that
\begin{equation}\label{inveq}
Q(a,x)=q
\end{equation}
holds, where $0<q<1$ and $a>0$ are given. The inversion of the normalised incomplete gamma function is an important topic in mathematical statistics and probability theory, in particular for computing percentage points of the $\chi^2$-distribution (see, for example, \cite{Fettis1979}). A suitable inversion of the expansion \eqref{eq7} yields the following result.

\begin{theorem}\label{thminv} For any $0<q<1$, define $\tau_0$ to be the unique real root of the equation
\begin{equation}\label{qdef}
q = \tfrac{1}{2}\erfc\left(2^{ - \frac{1}{2}} \tau _0\right).
\end{equation}
Then the inverse function $x=x(a,q)$ that satisfies the equation \eqref{inveq} has the asymptotic expansion
\begin{equation}\label{invexp}
x(a,q) \sim a + \tau _0 a^{\frac{1}{2}} + \sum_{n = 0}^\infty \frac{d_n(\tau_0)}{a^{n/2}},
\end{equation}
as $a\to +\infty$, provided that $q$ is bounded away from $0$ and $1$. The coefficients $d_n(\tau_0)$ are polynomials in $\tau_0$ of degree $n+2$. In addition, the even- and odd-order polynomials are even and odd functions, respectively.
\end{theorem}

The first few coefficients in \eqref{invexp} are as follows:
\[
d_0(\tau_0) = \frac{1}{3}\tau_{0}^2-\frac{1}{3},\quad d_1(\tau_0) = \frac{1}{36}\tau_{0}^3-\frac{7}{36}\tau_{0},\quad d_2(\tau_0) = -\frac{1}{270}\tau_{0}^4-\frac{7}{810}\tau_{0}^2+\frac{8}{405}.
\]
For higher coefficients and for an algorithm to generate them, see Section \ref{section4}.

In Theorem \ref{thminv}, we assume that $q$ is bounded away from $0$ and $1$. In practice, these bounds can be very small. It follows from \eqref{qdef} that $\tau_0 \to +\infty$ as $q\to 0$ and $\tau_0 \to -\infty$ as $q \to 1$. Since 
\[
\tfrac{1}{2}\erfc\left(2^{-\frac{1}{2}}\tau_0\right) \sim \frac{1}{\sqrt {2\pi }\tau_0}\exp\left(-\frac{\tau_0^2}{2}\right)
\]
as $\tau_0 \to +\infty$, and
\[
\tfrac{1}{2}\erfc\left(2^{-\frac{1}{2}}\tau_0\right) \sim 1+\frac{1}{\sqrt {2\pi }\tau_0}\exp\left(-\frac{\tau_0^2}{2}\right)
\]
as $\tau_0 \to -\infty$, we can infer that $|\tau _0|  \sim \sqrt { - 2\log (q(1 - q))} $ as $q\to 0$ or $1$. Since the $d_n(\tau_0)$'s are polynomials 
in $\tau_0$ of degree $n+2$, for \eqref{invexp} to be an asymptotic expansion we need $\tau_0=o\left(\left|a\right|^{1/2}\right)$ as $a\to\infty$. 
Hence, we need $\log (q(1 - q)) = o(|a|)$, i.e., $q$ can be sub-exponentially close in $\left|a\right|$ to $0$ or $1$.

The other expansion \eqref{eq8} is particularly well suited for the study of the large-$a$ behaviour of the $x$-zeros of the regularised 
incomplete gamma function $\gamma^\ast (a,x) = x^{-a} P(a,x) = x^{-a} (1 - Q(a,x))$. The function $\gamma^\ast(a,x)$ is entire in $a$ and $x$ 
and has a unique negative zero $x_{-}(a)$ when $a$ is negative and $a \neq 0,-1,-2,\ldots$. 
An asymptotic approximation for $x_{-}(a)$ as $a \to -\infty$ was given by Tricomi \cite{Tricomi1950}, which was corrected later 
by K\"{o}lbig \cite{Kolbig1972}. A more accurate asymptotic approximation was proved recently by Thompson \cite{Thompson2012}, which reads
\begin{equation}\label{eq11}
x_{-}(a) = a - \tau_1 \left(-a\right)^{\frac{1}{2}} - \tfrac{1}{3}\tau_1^2-\tfrac{1}{3} + \mathcal{O}\left(\left(-a\right)^{-\frac12}\right),
\end{equation}
as $a \to -\infty$, provided that $a$ is bounded away from the non-positive integers. Here $\tau_1$ is the unique solution of the equation
\begin{equation}\label{eq12}
\cot (- \pi a) = \sqrt {\frac{2}{\pi}} \int_0^{\tau _1 } \exp \left( \frac{t^2 }{2} \right)\,\d t  
=\frac{2}{\sqrt\pi}F\left(2^{-\frac{1}{2}}\tau_1\right)\exp\left(\frac{\tau_1^2}{2} \right) ,
\end{equation}
where $F(w)$ denotes Dawson's integral \cite[\href{http://dlmf.nist.gov/7.2.ii}{\S 7.2(ii)}]{NIST:DLMF}.
In the present paper, we extend the result \eqref{eq11} by deriving a complete asymptotic expansion for $x_-(a)$.

\begin{theorem}\label{thmzero} The unique negative zero $x_{-}(a)$ of the regularised incomplete gamma function $\gamma^\ast(a,x)$ has the asymptotic expansion
\begin{equation}\label{eq20}
x_-(a) \sim a - \tau_1 \left(-a\right)^{\frac{1}{2}} + \sum_{n = 0}^\infty \frac{(-\im)^n d_n(\im\tau_1)}{\left(-a\right)^{n/2}},
\end{equation}
as $a\to -\infty$, provided that $a$ is bounded away from the non-positive integers. Here $\tau_1$ is the unique solution of the equation 
\eqref{eq12}. The coefficients $d_n(\im\tau_1)$ are identical to those appearing in Theorem \ref{thminv} with $\tau_0$ replaced by 
$\im\tau_1$.
\end{theorem}

In Theorem \ref{thmzero}, we require that $a$ is bounded away from the non-positive integers. In practice, these bounds can be very small. Indeed, let $k$ be a negative integer. It follows from \eqref{eq12} that $\tau_1 \to \pm\infty$ as $a\to k$. Since 
\[
\cot(-\pi a) \sim  -\frac{1}{\pi(a - k)},
\]
as $a \to k$, and
\[
\frac{2}{\sqrt \pi}F\left(2^{-\frac{1}{2}} \tau_1\right)\exp \left( \frac{\tau_1^2}{2} \right) \sim \sqrt {\frac{2}{\pi}} \frac{1}{\tau_1}\exp \left( \frac{\tau _1^2 }{2} \right),
\]
as $\tau_1 \to \pm\infty$, we can infer from \eqref{eq12} that $|\tau_1| \sim \sqrt{-2\log | a-k |}$ as $a \to k$. 
Since the $d_n(\im \tau_1)$'s are polynomials in $\tau_1$ of degree $n+2$, for \eqref{eq20} to be an asymptotic expansion we 
need $\tau_1=o\left(\left|a\right|^{1/2}\right)$ as $a\to-\infty$. Thus, we need $\log |a - k| = o(|a|)$,  i.e., $a$ can be sub-exponentially 
close in $\left|a\right|$ to a negative integer.

A different asymptotic expansion for $x_-(a)$ can be derived using the known uniform asymptotic expansion of the function $\gamma^\ast(a,x)$ \cite{Temme1996b} and the method described in \cite{Temme1992}. The resulting expansion would also hold when $a$ is arbitrarily close to a non-positive integer. However, the coefficients in this expansion would be complicated functions possessing removable singularities, whereas the expansion \eqref{eq20} does not suffer from this inconvenience, because its coefficients are polynomials.

The remaining part of the paper is structured as follows. Section \ref{section1} provides a historical overview on the existing asymptotic 
expansions of the (normalised) incomplete gamma functions and the relation between those expansions and ours. In Section \ref{section2}, 
we prove the asymptotic expansions given in Theorem \ref{thm1}. In Section \ref{section3}, we show that Proposition \ref{prop1} holds. 
In Section \ref{section4}, we prove the asymptotic expansions for the inverse and for the unique negative zero of the regularised incomplete 
gamma function stated in Theorems \ref{thminv} and \ref{thmzero}. Finally,  in Section \ref{section5}, we compare numerically our 
new transition region expansion \eqref{eq7} with one of Temme's uniform asymptotic expansions. We focus especially on the case that $\tau$ is unbounded, but will make the surprising observation that even in the case that $\tau$ is comparable with $a^{1/6}$, we still obtain reasonable approximations from \eqref{eq7}.


\section{Historical overview}\label{section1}

The incomplete gamma functions of the complex variables $a$ and $z$ are defined by the integrals
\[
\gamma (a,z) = \int_0^z t^{a - 1} \e^{ - t} \,\d t ,\quad \Gamma(a,z) = \int_z^\infty t^{a - 1} \e^{ - t} \,\d t ,
\]
where the paths of integration do not cross the negative real axis, and in the case of $\Gamma(a,z)$ exclude the origin. The definition of $\gamma(a,z)$ requires the condition $\Re(a) > 0$ while in that for $\Gamma(a,z)$ it is assumed that $|\arg z|<\pi$.  When $z\neq 0$, $\Gamma(a,z)$ is an entire function of $a$, and $\gamma(a,z)$ extends to a meromorphic function with simple poles at $a=-n$, $n \in \mathbb{N}$, with residue $(-1)^n/n!$. For fixed values of $a$, both $\gamma(a,z)$ and $\Gamma(a,z)$ are multivalued functions of $z$ and are holomorphic on the Riemann surface of the complex logarithm.

The asymptotic behaviour of $\gamma(a,z)$ and $\Gamma(a,z)$ when either $a$ or $z$ is large is well understood (see, for instance, 
\cite[\href{http://dlmf.nist.gov/8.11.i}{\S 8.11(i)} and \href{http://dlmf.nist.gov/8.11.ii}{\S 8.11(ii)}]{NIST:DLMF}).
The treatment when both $a$ and $z$ become large is significantly more complicated since the resulting expansions have to take into account the presence of the transition point at $z=a$, about which the asymptotic structures of $\gamma(a,z)$ and $\Gamma(a,z)$ go through an abrupt change. For example, in the case of large positive variables, the normalised incomplete gamma function $Q(a,z)=\Gamma(a,z)/\Gamma(a)$ exhibits a sharp decay near the transition point $z=a$, since it is approximately unity when $z<a$ and decreases algebraically to zero when $z>a$. The behaviour of the complementary normalised incomplete gamma function $P(a,z)=\gamma(a,z)/\Gamma(a)$ follows from the functional relation $P(a,z)+Q(a,z)=1$.

The earliest asymptotic expansions of $\gamma(a,z)$ and $\Gamma(a,z)$ for large $a$ and $z$ were given by Mahler \cite{Mahler1930} and Tricomi \cite{Tricomi1950}. In modern notation, taking $\lambda=z/a$, Mahler's expansion for $\gamma(a,z)$ may be written
\begin{equation}\label{eq1}
\gamma (a,z) \sim - z^a \e^{ - z} \sum_{n = 0}^\infty \frac{\left( - a\right)^n b_n (\lambda )}{\left(z - a\right)^{2n + 1} }, 
\end{equation}
as $a\to \infty$ in the sector $|\arg a|\leq \frac{\pi}2-\delta<\frac{\pi}2$, provided $0<\lambda<1$. With the same notation, Tricomi's result for $\Gamma(a,z)$ reads
\begin{equation}\label{eq2}
\Gamma (a,z) \sim z^a \e^{ - z} \sum_{n = 0}^\infty \frac{\left( - a\right)^n b_n (\lambda )}{\left(z - a\right)^{2n + 1} } ,
\end{equation}
as $a\to \infty$ in the sector $|\arg a|\leq \frac{3\pi}{2}-\delta<\frac{3\pi}{2}$, provided $\lambda>1$. The coefficients $b_n(\lambda)$ are polynomials in $\lambda$ of degree $n$, the first few being
\begin{gather*}
b_0 (\lambda ) = 1,\quad b_1 (\lambda ) = \lambda ,\quad b_2 (\lambda ) = 2\lambda ^2  + \lambda ,\\ b_3 (\lambda ) = 6\lambda ^3  + 8\lambda ^2  + \lambda ,\quad b_4 (\lambda ) = 24\lambda ^4  + 58\lambda ^3  + 22\lambda ^2  + \lambda .
\end{gather*}
Higher coefficients can be computed using the recurrence relation \cite[\href{http://dlmf.nist.gov/8.11.9}{Eq.~8.11.9}]{NIST:DLMF}
\[
b_k(\lambda)=\lambda(1-\lambda)b'_{k-1}(\lambda)+(2k-1)\lambda b_{k-1}(\lambda) ,
\]
for $k\geq 1$. For other representations of $b_k(\lambda)$, including an explicit expression involving the Stirling numbers of the second kind, the reader is referred to \cite{Nemes16}.

An analogous expansion for $\Gamma (-a,z)$ was subsequently provided by Gautschi \cite{Gautschi1959} (which was extended later to complex $a$ by Temme \cite{Temme1994}) in the form
\begin{equation}\label{eq3}
\Gamma (-a,z) \sim z^{-a} \e^{ - z} \sum_{n = 0}^\infty \frac{a^n b_n (-\lambda )}{\left(z + a\right)^{2n + 1} } ,
\end{equation}
as $a\to \infty$ in the sector $|\arg a|\leq \frac{\pi}{2}-\delta< \frac{\pi}{2}$, provided $z/a=\lambda>-1$. If we restrict $\lambda$ to positive values,
 the expansion \eqref{eq3} remains valid in the wider range $|\arg a| \le \frac{3\pi}{2} - \omega  - \delta < \frac{3\pi}{2} - \omega$ with 
 $\omega  = \arg (  \lambda  + \log \lambda + \pi\im)$, $0 < \omega  < \pi$. For large positive values of $a$ and complex values of $\lambda$, 
 the asymptotic behaviour of $\Gamma (-a,z)$ was discussed by Dunster \cite{Dunster1996}.

The asymptotic expansions \eqref{eq1}--\eqref{eq2} break down as $\lambda \to 1$ (i.e., as the transition point $z=a$ is approached) since 
their terms become singular in this limit. Similarly, \eqref{eq3} fails to hold as $\lambda \to -1$. An expansion which is useful near the transition 
point is as follows:
\begin{equation}\label{eq4}
\Gamma (a,z) \sim z^a \e^{ - z} \sqrt {\frac{\pi }{2z}} \sum_{n = 0}^\infty  \frac{a_{2n} (\varepsilon )}{z^n }  - z^{a - 1} \e^{ - z} \sum_{n = 0}^\infty \frac{a_{2n + 1} (\varepsilon )}{z^n } ,
\end{equation}
as $z\to \infty$ in the sector $|\arg z|\leq 2\pi-\delta<2\pi$ with $\varepsilon = z - a$ being bounded. The coefficients $a_n(\varepsilon)$ are polynomials in $\varepsilon$ of degree $n$, the first few of which are
\begin{gather*}
a_0 (\varepsilon ) = 1,\quad a_1 (\varepsilon ) =  \varepsilon  + \frac{1}{3},\quad a_2 (\varepsilon ) = \frac{1}{2}\varepsilon ^2  + \frac{1}{2}\varepsilon  + \frac{1}{{12}},\\  a_3 (\varepsilon ) = \frac{1}{3}\varepsilon ^3  + \frac{2}{3}\varepsilon ^2  + \frac{1}{3}\varepsilon  + \frac{4}{{135}},\quad a_4 (\varepsilon ) = \frac{1}{8}\varepsilon ^4  + \frac{5}{{12}}\varepsilon ^3  + \frac{5}{{12}}\varepsilon ^2  + \frac{1}{8}\varepsilon  + \frac{1}{{288}}.
\end{gather*}
In the special case that $\varepsilon=0$, the expansion \eqref{eq4} is well known. Its remarkably large region of validity was proved recently by the first author \cite{Nemes16}. For the case of general complex $\varepsilon$, \eqref{eq4} can be proved by employing the method of steepest descents to the integral representation \cite[\href{http://dlmf.nist.gov/8.6.7}{Eq.~8.6.7}]{NIST:DLMF}
\[
\Gamma (a,z) = z^a \e^{ - z} \int_0^{ + \infty } \exp \left( - z(\e^t  - t - 1)\right)\e^{ - \varepsilon t} \,\d t, \quad \Re(z)>0.
\]
The region of validity of the resulting asymptotic expansion \eqref{eq4} depends only on the singularity structure of the inverse function of $\e^t  - t - 1$ and hence, it is independent of $\varepsilon$. 

Dingle \cite[Ch. VIII, Sec. 3 and Ch. XXIII, Sec. 8]{Dingle73} gave several interesting (formal) results about the expansions \eqref{eq1}--\eqref{eq3} and \eqref{eq4} with $\varepsilon=0$, including re-expansions for the remainder terms and asymptotic approximations for the high-order coefficients. Many of his results have been rigorously justified recently by the first author \cite{Nemes15b,Nemes16}. These papers also contain sharp bounds for the error terms of the asymptotic expansions \eqref{eq2}, \eqref{eq3} and \eqref{eq4} with $\varepsilon=0$.

The expansions \eqref{eq1}--\eqref{eq3} are suitable for numerical computation when $|z-a|$ is large compared with $|a|^{1/2}$, for otherwise 
the early terms of these expansions do not decrease in magnitude. 
Similarly, the expansion \eqref{eq4} is useful only when $z = a + o\left(\left|a\right|^{1/2}\right)$. Thus, corresponding to values of $|z-a|$ that are 
comparable with $\left|a\right|^{1/2}$, there are gaps within which neither of these expansions is suitable. We call these gaps the transition regions.

The first attempt to obtain an asymptotic estimate which is valid in the transition regions was made by Tricomi \cite{Tricomi1950}, who showed that
\[
Q\left(a + 1,a + \tau a^{\frac12} \right) = \tfrac{1}{2}\erfc\left(2^{ - \frac{1}{2}} \tau \right) + 
\frac{1}{\sqrt {2\pi a}}\exp \left( -\frac{\tau ^2 }{2} \right)\frac{\tau ^2  + 2}{3} + \mathcal{O}\left(\left|a\right|^{-\frac12} \right),
\]
as $a\to \infty$ in the sector $|\arg a|\leq \frac{\pi}{2}-\delta< \frac{\pi}{2}$, with any fixed $\tau$ satisfying 
$\left|\arg \left(\tau a^{1/2} \right)\right| \le \pi  - \delta  < \pi$. He also gave, in place of the error term $\mathcal{O}\left(\left|a\right|^{-1/2}\right)$, 
a complicated expansion in terms of incomplete gamma functions.

A similar asymptotic approximation for $Q\left(a,a + \tau a^{1/2}\right)$ was subsequently given by Pagurova \cite{Pagurova1956} with an error 
term of $\mathcal{O}\left(a^{-3}\right)$. Pagurova's approximation is limited to positive real values of $a$ and real values of $\tau$, but is considerably 
simpler than that obtained by Tricomi.

Expansions of a different type were provided by Paris \cite{Paris2002b,Paris2018} in the form
\begin{equation}\label{eq4a}
\left.\begin{array}{l}\gamma(a,z)\\\Gamma(a,z)\end{array}\!\!\right\} \sim z^{a - \frac{1}{2}} \e^{ -z} \left( \sqrt {\frac{\pi}{2}} 
\exp \left( \frac{\chi ^2 }{2} \right)\erfc \left( \mp 2^{ - \frac{1}{2}} \chi \right)\sum_{n = 0}^\infty \frac{A_n (\chi)}{z^{n/2}}  \pm 
\sum_{n = 1}^\infty \frac{B_n (\chi )}{z^{n/2}}  \right)
\end{equation}
as $z\to \infty$ in the sector $|\arg z|\leq \frac{\pi}{2}-\delta< \frac{\pi}{2}$, and $\Re(z - a) \leq 0$ for $\gamma(a, z)$ and $\Re(z - a) \geq 0$ for $\Gamma(a, z)$. Here $\chi$ is defined by the equality $a = z - \chi z^{1/2}$ and the coefficients $A_n(\chi)$ and $B_n(\chi)$ are polynomials in $\chi$ of degree $3n$ and $3n-1$, respectively. The first few of these coefficients are as follows:
\begin{gather*}
A_0 (\chi ) = 1,\quad A_1 (\chi ) = \frac{1}{6}\chi ^3  + \frac{1}{2}\chi ,\quad A_2 (\chi ) = \frac{1}{72}\chi ^6  + \frac{1}{6}\chi ^4  + \frac{3}{8}\chi ^2  + \frac{1}{12},\\
B_1 (\chi ) = \frac{1}{6}\chi ^2  + \frac{1}{3},\quad B_2 (\chi ) = \frac{1}{72}\chi ^5  + \frac{11}{72}\chi ^3  + \frac{1}{4}\chi .
\end{gather*}
The asymptotic expansions \eqref{eq4a} together cover the transition regions and are also valid away from those regions, i.e., when $\chi$ is large. However, as it was noted by Paris \cite{Paris2018}, even for moderately large values of $\chi$, the practical use of the expansions is problematic because of severe numerical cancellations. Therefore, the use of these asymptotic expansions should be confined to bounded values of $\chi$.

Other asymptotic expansions for $\gamma(a,z)$ and $\Gamma(a,z)$, similar to those of \eqref{eq4a}, were derived by Dingle \cite[p. 249]{Dingle73} and L\'opez et al. \cite{FLP05}.

Asymptotic expansions for the normalised incomplete gamma function $Q(a,z)$ which are uniformly valid in the variables $a$ and $z$ were first established by Temme \cite{Temme1979,Temme1996b}. His results may be stated as follows. Define
\[
\lambda  = z/a \; \text{ and } \; \eta =\eta(\lambda) = \left(2(\lambda  - 1 - \log \lambda )\right)^{\frac{1}{2}} ,
\]
where the branch of the square root is continuous and satisfies $\eta(\lambda) \sim \lambda-1$ as $\lambda \to 1$. Then as $a\to \infty$ in the sector $|\arg a|\leq \pi-\delta< \pi$,
\begin{equation}\label{eq5}
Q(a,z) \sim \tfrac{1}{2} \erfc\left(2^{ - \frac{1}{2}} \eta a^{\frac{1}{2}} \right) + \frac{1}{\sqrt {2\pi a} }\exp \left( -\tfrac{1}{2}\eta^2 a \right)
\sum_{n = 0}^\infty \frac{c_n(\eta)}{a^n}
\end{equation}
and
\begin{equation}\label{eq6}
\frac{\e^{ \pm \pi \im a}}{2\im \sin (\pi a)}Q\left( - a,z\e^{ \pm\pi\im}\right) \sim  \pm \tfrac{1}{2} \erfc\left(\pm 2^{-\frac{1}{2}}\im\eta a^{\frac{1}{2}} \right) 
- \frac{\im }{\sqrt {2\pi a}}\exp \left( {\tfrac{1}{2}\eta ^2 a} \right)\sum_{n = 0}^\infty \left(-1\right)^n \frac{c_n (\eta )}{a^n } ,
\end{equation}
uniformly with respect to $\lambda$ in the sector $|\arg \lambda|\leq 2\pi-\delta< 2\pi$. The coefficients $c_n (\eta)$ are holomorphic functions of $\eta$ and are defined recursively by
\begin{equation}\label{eq6a}
c_0 (\eta ) = \frac{1}{\lambda  - 1} - \frac{1}{\eta},\quad c_n (\eta) = \frac{\gamma _n }{\lambda  - 1} + \frac{1}{\eta}\frac{\d c_{n - 1} (\eta )}{\d\eta}\quad (n \ge 1),
\end{equation}
where the $\gamma_n$'s are the Stirling coefficients appearing in the well-known asymptotic expansion of the gamma function (see, for example, \cite{Nemes13}). 
Temme's expansions feature the characteristic error function behaviour near the transition point $z=a$ (where $\eta=0$) and also describe the large-$a$ asymptotics uniformly in $z$. 
However, they suffer from the inconvenience of the coefficients $c_n (\eta)$ having a removable singularity at $\eta=0$, which makes their evaluation in the neighbourhood of the transition point difficult. 
To overcome this difficulty, Temme \cite{Temme1979} gave power series expansions in $\eta$ for these coefficients which do not have a removable singularity. 
The asymptotic properties of the $c_n (\eta)$'s for large $n$ were studied by Dunster et al. \cite{Dunster1998} and the second author \cite{OD98a}. 
Explicit expressions for these coefficients can be found in the paper \cite{Nemes13}. 
Computable error bounds for the asymptotic expansion \eqref{eq5} were established by Temme \cite{Temme1979} for real variables, and by Paris \cite{Paris2002a} for complex variables.

Although our expansions \eqref{eq7} and \eqref{eq8} are not valid in as large a domain of $a$ and $z$ as those in \eqref{eq5} and \eqref{eq6}, the coefficients have the advantage of not possessing a removable singularity at the transition point $z=a$ and so are more straightforward to compute. Pagurova's result is a special case of the expansion \eqref{eq7} truncated after the first six terms and restricted to large positive real values of $a$ and real values of $\tau$. The expansions \eqref{eq4a} are similar in character to the expansion \eqref{eq7}, but they are for the non-normalised functions $\gamma(a,z)$ and $\Gamma(a,z)$ and have more complicated forms and smaller regions of validity compared to \eqref{eq7}. Also, the coefficients of these asymptotic expansions seem not to satisfy simple recurrence relations like those of the expansion \eqref{eq7}. Asymptotic expansions analogous to \eqref{eq7} and \eqref{eq8} for the complementary normalised incomplete gamma function $P(a,z)$ follow from the functional relation $P(a,z)+Q(a,z)=1$.


\section{Proof of Theorem \ref{thm1}}\label{section2}

From \eqref{eq5}, we can assert that for any non-negative integer $N$,
\[
Q(a,z) = \tfrac{1}{2}\erfc\left(2^{ - \frac{1}{2}} \eta a^{\frac{1}{2}}\right)+ \frac{1}{\sqrt {2\pi a}}\exp \left(  - \tfrac{1}{2}\eta ^2 a \right)\left( 
\sum_{n = 0}^{N - 1} \frac{c_n (\eta )}{a^n }  + \mathcal{O}_{N,\delta } \left( \frac{1}{|a|^N } \right) \right),
\]
as $a\to \infty$ in the sector $|\arg a|\leq \pi-\delta< \pi$, uniformly with respect to $\lambda$ in the sector $|\arg \lambda|\leq 2\pi-\delta< 2\pi$. 
(Throughout this paper, we use subscripts in the $\mathcal{O}$ notations to indicate the dependence of the implied constant on certain parameters.) 
We employ this result with the choice of $\lambda  = 1 + \tau a^{ - 1/2}$ and assume that $\left|\tau a^{ - 1/2}\right|<\frac{1}{2}$, 
in particular, $|\arg \lambda|<\frac{\pi}{2}$. Consequently, we have
\begin{multline*}
\tfrac{1}{2}\left|\eta\right|^2  = \left| \tau a^{ - 1/2}  - \log (1 + \tau a^{ - 1/2} )\right| 
=\left| \tau a^{ - 1/2}\right|^2 \left| \int_0^1 \frac{t}{1 + \tau a^{ - 1/2} t}\,\d t \right| \\
< \left| \tau a^{ - 1/2} \right|^2 \int_0^1 \frac{t}{1 - \left| {\tau a^{ - 1/2} } \right|t}\,\d t  < \left| \tau a^{ - 1/2} \right|^2  < \left| \tau a^{ - 1/2} \right|,
\end{multline*}
that is, $|\eta| < 1$. We have the convergent expansions
\begin{align*}
\tfrac{1}{2}\eta ^2  = \tau a^{ - 1/2}  - \log \left(1 + \tau a^{ - 1/2} \right) &= \sum_{k = 2}^\infty\left(-1\right)^k \frac{\tau ^k}{k}\frac{1}{a^{k/2}} 
 \\& = \frac{\tau ^2 }{2}\frac{1}{a}\left( 1 + \sum_{k = 1}^\infty\left(-1\right)^k \frac{2\tau ^k }{k + 2}\frac{1}{a^{k/2} } \right)
\end{align*}
and
\[
\eta  = \frac{\tau }{a^{1/2}} + \sum_{k = 2}^\infty \frac{m_k (\tau )}{a^{k/2}},
\]
where $m_k(\tau)$ is a monomial in $\tau$ of degree $k$.

Let $h$ denote a complex number whose value will be specified later. Employing the known expression for the higher derivatives of the
complementary error function \cite[\href{http://dlmf.nist.gov/7.10.1}{Eq.~7.10.1}]{NIST:DLMF} and the Taylor formula with integral remainder, we find
\begin{align*}
\tfrac{1}{2}\erfc\left(2^{ - \frac{1}{2}} \tau  + h\right)=\; & \tfrac{1}{2}\erfc\left(2^{ - \frac{1}{2}} \tau\right)\\&+ \frac{1}{\sqrt \pi}\exp \left( - \frac{\tau ^2}{2}\right)\left( \sum_{n = 1}^{N - 1} \frac{\left(-1\right)^n}{n!}H_{n - 1} \left(2^{ - \frac{1}{2}} \tau \right)h^n   + R_N (h,\tau ) \right),
\end{align*}
where
\begin{align*}
R_N (h,\tau ) &= \frac{\sqrt \pi}{2(N - 1)!}\exp \left( \frac{\tau ^2}{2} \right)\int_0^h {\erfc^{(N)}}\left(2^{ - \frac{1}{2}} \tau+t\right)\left(h-t\right)^{N-1} \,\d t 
\\ &= \frac{\left(-1\right)^N }{(N - 1)!}\int_0^h H_{N - 1}\left(2^{ - \frac{1}{2}} \tau+t\right)\exp\left(-2^{\frac{1}{2}} \tau t-t^2 \right)\left(h-t\right)^{N-1} \,\d t
\\ &= \mathcal{O}_N\left(\left(|\tau|  + 1\right)^{N - 1} \left|h\right|^N \exp \left(\mathcal{O}(|\tau h|)\right) \right)
\end{align*}
as $h\to 0$ and $H_n(w)$ denotes the $n$th Hermite polynomial \cite[\href{http://dlmf.nist.gov/18.3}{\S 18.3}]{NIST:DLMF}.
Assuming that $\tau  = o\left(|a|^{1/4} \right)$ for large $a$ and using this expansion with
\[
h = \sum_{k = 1}^\infty  2^{ - \frac{1}{2}} \frac{m_{k + 1} (\tau )}{a^{k/2} }  = \mathcal{O}\left( \frac{|\tau |^2 }{|a|^{1/2} } \right) = o(1)
\]
and with $N+1$ in place of $N$, we obtain
\begin{align*}
& \tfrac{1}{2}\erfc\left(2^{ - \frac{1}{2}} \eta a^{\frac{1}{2}}\right)= \tfrac{1}{2}\erfc \left( 2^{ - \frac{1}{2}} \tau  + \sum_{k = 1}^\infty  2^{ - \frac{1}{2}} \frac{m_{k + 1} (\tau )}{a^{k/2}}  \right)  = \tfrac{1}{2}\erfc\left(2^{ - \frac{1}{2}} \tau\right)\\ &+ \frac{1}{\sqrt {2\pi a}}\exp \left( { - \frac{\tau ^2 }{2}} \right)\left( 
\sum_{n = 0}^{N - 1} \frac{p_n (\tau )}{a^{n/2}}  + \mathcal{O}_N \left( \frac{(|\tau|  + 1)^{3N + 2}}{\left| a \right|^{N/2} }
\exp \left( \mathcal{O}\left( \frac{|\tau |^3 }{|a|^{1/2} } \right) \right) \right) \right).
\end{align*}
Here $p_n(\tau)$ is a polynomial in $\tau$ of degree at most $3n+2$.

We have
\[
\exp \left(  - \frac{\tau ^2}{2} + h \right) = \exp \left( - \frac{\tau ^2}{2} \right)\left( {\sum_{n = 0}^{N - 1} {\frac{{h^n }}{{n!}}}  
+ \mathcal{O}_N \left(|h|^N \exp (|h|)\right)} \right)
\]
for any complex $h$. Applying this expansion with
\[
h = \sum_{k = 1}^\infty\left(-1\right)^{k + 1} \frac{\tau ^{k + 2} }{k + 2}\frac{1}{a^{k/2}}  = 
\mathcal{O}\left( \frac{\left|\tau\right|^3}{\left|a\right|^{1/2}} \right),
\]
we deduce
\begin{align*}
\exp \left(  - \tfrac{1}{2}\eta ^2 a \right) & = \exp \left(  - \frac{\tau ^2}{2} + \sum_{k = 1}^\infty\left(-1\right)^{k + 1} 
\frac{\tau ^{k + 2}}{k + 2}\frac{1}{a^{k/2}} \right) \\ 
& = \exp \left(  - \frac{\tau ^2}{2} \right)\left( \sum_{n = 0}^{N - 1} \frac{q_n (\tau )}{a^{n/2}}  
+ \mathcal{O}_N \left( {\frac{\left(|\tau |+1\right)^{3N} }{\left|a\right|^{N/2} }
\exp \left( \mathcal{O}\left( \frac{\left|\tau\right|^3 }{\left|a\right|^{1/2} } \right) \right)} \right) \right),
\end{align*}
where $q_n(\tau)$ is a polynomial in $\tau$ of degree at most $3n$.

We have the power series expansion
\[
c_n (h) = c_n (0) + \sum_{k = 1}^\infty  {f_{n,k} h^k } ,
\]
which converges if $|h|<2\sqrt{\pi}$ (cf. \cite[\href{http://dlmf.nist.gov/8.12.11}{Eq.~8.12.11}]{NIST:DLMF} or Appendix \ref{appendixb}). 
Applying this expansion with $h= \eta$ ($|\eta|<1$), we find
\[
c_n (\eta ) = c_n \left( \frac{\tau }{a^{1/2}} + \sum_{k = 2}^\infty \frac{m_k (\tau )}{a^{k/2} } \right) = \sum_{k = 0}^\infty \frac{r_k (\tau )}{a^{k/2} } ,
\]
where $r_k(\tau)$ is a monomial in $\tau$ of degree $k$. Consequently,
\[
\sum_{n = 0}^{N - 1} \frac{c_n (\eta )}{a^n}  + \mathcal{O}_{N,\delta } \left( \frac{1}{|a|^N }\right) = 
\sum_{n = 0}^{N - 1} \frac{s_n (\tau )}{a^{n/2}}  + \mathcal{O}_{N,\delta } \left( \frac{\left(|\tau | + 1\right)^N }{\left|a\right|^{N/2} } \right)
\]
where $s_n(\tau)$ is a polynomial in $\tau$ of degree at most $n$. Collecting all the partial results, we finally have
\begin{gather}\label{eq22}
\begin{split}
& Q(a,z)=Q\left(a,a+\tau a^{\frac{1}{2}}\right) = \tfrac{1}{2} \erfc\left(2^{-\frac{1}{2}}\tau\right) \\&+ \frac{1}{\sqrt {2\pi a}}\exp\left(-\frac{\tau^2}{2}\right)
\left(\sum_{n=0}^{N-1}\frac{C_n(\tau)}{a^{n/2}} + \mathcal{O}_{N,\delta} \left( \frac{\left(|\tau|  + 1\right)^{3N + 2}}{\left| a \right|^{N/2} }
\exp \left( {\mathcal{O}\left( \frac{\left|\tau\right|^3 }{\left|a\right|^{1/2} } \right)} \right) \right) \right)
\end{split}
\end{gather}
for any non-negative integer $N$ as $a\to \infty$ in the sector $|\arg a|\leq \pi-\delta< \pi$, provided that $\tau  = o\left(\left|a\right|^{1/4}\right)$. 
Here $C_n(\tau)$ is a polynomial in $\tau$ of degree at most $3n+2$. If $\tau$ is bounded, then the error term in \eqref{eq22} is of the same 
order of magnitude as the first neglected term and hence, the asymptotic expansion \eqref{eq7} holds.

In a similar manner, starting with the asymptotic expansion \eqref{eq6}, it can be shown that
\begin{gather}\label{eq9}
\begin{split}
& \frac{\e^{\pm\pi\im a}}{2\im\sin (\pi a)}Q\left(-a,(a+\tau a^{\frac{1}{2}})\e^{\pm\pi\im}\right) = \pm \tfrac{1}{2}\erfc\left(\pm2^{-\frac{1}{2}}\im\tau\right) \\ 
& - \frac{\im}{\sqrt{2\pi a}}\exp \left(\frac{\tau^2}{2}\right) \left(\sum_{n=0}^{N-1}\frac{\widetilde{C}_n(\tau)}{a^{n/2}} 
+ \mathcal{O}_{N,\delta} \left( \frac{\left(|\tau|  + 1\right)^{3N + 2}}{\left| a \right|^{N/2} }
\exp \left( {\mathcal{O}\left( \frac{\left|\tau\right|^3 }{\left|a\right|^{1/2} } \right)} \right) \right) \right)
\end{split}
\end{gather}
for any non-negative integer $N$ as $a\to \infty$ in the sector $|\arg a|\leq \pi-\delta< \pi$, provided that $\tau  = o\left(\left|a\right|^{1/4}\right)$. 
The coefficients $\widetilde{C}_n(\tau)$ are polynomials in $\tau$ of degree at most $3n+2$. If $\tau$ is bounded, then the error term in \eqref{eq9} 
is of the same order of magnitude as the first neglected term and hence, the asymptotic expansion \eqref{eq8} holds. 
To show that $\widetilde C_n (\tau ) = \left(-\im\right)^n C_n (\im\tau )$, we can proceed as follows. Assuming that $\tau$ is bounded and applying 
\eqref{eq9} with $a\e^{ \mp \pi\im/2}$ in place of $a$, $a \to +\infty$, we find
\begin{multline*}
\frac{1}{1 - \e^{ - 2\pi a}}Q\left( a\e^{ \pm \frac{\pi }{2}\im} ,a\e^{ \pm \frac{\pi}{2}\im}  + 
\tau \e^{ \pm \frac{\pi}{2}\im} \left( a\e^{ \pm \frac{\pi}{2}\im}  \right)^{\frac{1}{2}}  \right) \\ 
\sim \tfrac{1}{2} \erfc\left(2^{ - \frac{1}{2}} \tau \e^{ \pm \frac{\pi}{2}\im}\right) + \frac{1}{\sqrt {2\pi a\e^{ \pm \frac{\pi }{2}\im} } }
\exp \left( \frac{\tau ^2 }{2} \right)\sum_{n = 0}^\infty \frac{\left( \pm\im\right)^n \widetilde{C}_n (\tau)}{\left(a\e^{ \pm \frac{\pi}{2}\im}\right)^{n/2}} .
\end{multline*}
Since, for any non-negative integer $N$, $1/(1 - \e^{ - 2\pi a} ) = 1+\mathcal{O}(a^{ - N} )$ as $a\to +\infty$, the right hand side must agree with that 
of \eqref{eq5} when applied with $\tau \e^{ \pm \frac{\pi }{2}\im}$ in place of $\tau$ and $a\e^{ \pm \frac{\pi }{2}\im}$ in place of $a$, $a \to +\infty$. 
Therefore, we find that $\left( \pm\im\right)^n \widetilde C_n (\tau ) = C_n ( \pm\im\tau )$, that is 
$\widetilde C_n (\tau ) = \left( \mp\im\right)^n C_n ( \pm \im\tau )$. 
By the parity properties of the polynomials $C_n(\tau)$ (see below), this latter equality simplifies to 
$\widetilde C_n(\tau)=\left(-\im\right)^nC_n(\im\tau )$.

We proceed by proving the claimed properties of the polynomials $C_n(\tau)$. The differential equation 
\cite[\href{http://dlmf.nist.gov/8.2.iii}{\S 8.2(iii)}]{NIST:DLMF}
\[
\frac{\partial ^2 Q(a,z)}{\partial z^2} + \left( 1 + \frac{1 - a}{z} \right)\frac{\partial Q(a,z)}{\partial z} = 0
\]
implies that
\[
\left(1 + \tau a^{ - \frac{1}{2}}\right)\frac{\partial ^2 Q\left(a,a + \tau a^{\frac{1}{2}}\right)}{\partial \tau ^2 } 
+ \left(\tau  + a^{ - \frac{1}{2}} \right)\frac{\partial Q\left(a,a + \tau a^{\frac{1}{2}} \right)}{\partial \tau } = 0.
\]
Substituting \eqref{eq7} into the left-hand side of this equation and equating the coefficients of $a^{-n/2}$ on both sides yield $C_0 (\tau ) = \frac{1}{3}\tau^2 -\frac{1}{3}$ and
\begin{equation}\label{eq10}
C_n (\tau ) + \tau C'_n (\tau ) - C''_n (\tau ) = \tau (\tau ^2  - 2)C_{n - 1} (\tau ) - (2\tau ^2  - 1)C'_{n - 1} (\tau ) + \tau C''_{n - 1} (\tau )
\end{equation}
for $n\geq 1$. Using induction on $n$ and the fact that the corresponding homogeneous equation $w(\tau ) + \tau w'(\tau ) - w''(\tau ) = 0$ has no non-zero polynomial solution, it follows that \eqref{eq10} determines the polynomials $C_n(\tau)$ uniquely. An induction on $n$ and \eqref{eq10} imply that the degree of the polynomial $C_n(\tau)$ is exactly $3n+2$.

Lastly, we prove that $C_n ( - \tau )= \left(-1\right)^nC_n ( \tau )$, i.e., the even- and odd-order polynomials are even and odd functions, respectively. We substitute the polynomial expansion
\begin{equation}\label{eq10a}
C_n(\tau) = \sum_{k=0}^{3n+2} c_{n,k}\tau^k,
\end{equation}
into \eqref{eq10} and obtain for the coefficients the recurrence relation
\begin{equation}\label{eq10b}
c_{n,k}=(k+2)c_{n,k+2}+(k+1) c_{n-1,k+1}-\frac{2k}{k+1}c_{n-1,k-1}+\frac{1}{k+1}c_{n-1,k-3}.
\end{equation}
Since $C_0 (\tau ) = \frac{1}{3}\tau^2 -\frac{1}{3}$, it follows that
\begin{equation}\label{eq10c}
c_{n,3n+2}=\frac{1}{3^{n+1}(n+1)!},\qquad c_{n,3n+1}=0,
\end{equation}
holds for $n=0$, and by the recurrence relation \eqref{eq10b} we can show that \eqref{eq10c} holds for $n\geq1$. Starting with the first equation in \eqref{eq10c} and then taking $k=3n, 3n-2, 3n-4, \ldots$, we can compute the non-zero coefficients in the expansion \eqref{eq10a}. Similarly, it follows that $c_{n,k}=0$ for $k=3n+1, 3n-1, 3n-3, \ldots$. Accordingly, $C_n ( - \tau )= \left(-1\right)^nC_n ( \tau )$.

The first several polynomials $C_n ( \tau )$ are given in Table \ref{table1}.

\begin{table*}[!ht]
\begin{center}
\begin{tabular}
[c]{l l}  \hline
& \\ [-2ex]
 $n$ & $C_n(\tau)$ \\ [-2ex]
& \\ \hline 
 & \\ [-1.5ex]
 0 & $\frac{1}{3}{\tau}^{2}-\frac{1}{3}$\\ [1ex]
 1 & $\frac{1}{18}{\tau}^{5}-{\frac {11}{36}} {\tau}^{3}+\frac{1}{12}\tau$\\ [1ex]
 2 & $\frac{1}{162}\tau^{8}-\frac{29}{324}\tau^{6}+\frac{133}{540}\tau^{4}-\frac{23}{540}\tau^{2}-\frac{1}{540}$\\[1ex]
 3 & ${\frac {1}{1944}}{\tau}^{11}-{\frac {7}{486}}{\tau}^{9}+{
\frac {463}{4320}}{\tau}^{7}-{\frac {883}{4320}}{\tau}^{5}+{\frac 
{23}{864}}{\tau}^{3}-{\frac {1}{288}}\tau$\\[1ex]
 4 & ${\frac {1}{29160}}{\tau}^{14}-{\frac {23}{14580}}{
\tau}^{12}+{\frac {881}{38880}}{\tau}^{10}-{\frac {1507}{12960}}{
\tau}^{8}+{\frac {7901}{45360}}{\tau}^{6}-{\frac {61}{3024}}{\tau
}^{4}+{\frac {23}{6048}}{\tau}^{2}+{\frac{25}{6048}}$\\[1ex]
 5 & ${\frac {1}{524880}}{\tau}^{17}-{\frac {137}{
1049760}}{\tau}^{15}+{\frac {2149}{699840}}{\tau}^{13}-{\frac {42331}{1399680}}{\tau}^{
11}+{\frac {5910101}{48988800}}{\tau}^{9}-{\frac {413177
}{2721600}}{\tau}^{7}$\\[1ex]
& $+{\frac {3239}{194400}}{\tau}^{5}-{\frac {319
}{155520}}{\tau}^{3}-{\frac {139}{51840}}\tau$ \\[1ex]
 6 & ${\frac {1}{11022480}}{\tau}^{20}-{\frac {
191}{22044960}}{\tau}^{18}+{\frac {1483}{4898880}}{\tau}^{16}-{
\frac {47639}{9797760}}{\tau}^{14}+{\frac {1810717}{
48988800}}{\tau}^{12}$\\[1ex]
 & $-{\frac {1996859}{16329600}}{\tau}^{10}+{\frac {146611}{1088640}}{
\tau}^{8}-{\frac {2203}{155520}}{\tau}^{6}+{\frac {11}{10368}}{
\tau}^{4}+{\frac {259}{155520}}{\tau}^{2}+{\frac{101}{155520}}$ \\[1ex]
 7 & ${\frac {1}{264539520}}{\tau}^{23}-{\frac {127}{
264539520}}{\tau}^{21}+{\frac {2939}{125971200}}{\tau}^{19}-{\frac {2911}{5248800}}{\tau
}^{17}+{\frac {3217793}{470292480}}{\tau}^{15}$\\[1ex]
& $-{\frac {
33561391}{783820800}}{\tau}^{13}+{\frac {13697119}{
111974400}}{\tau}^{11}-{\frac {3797293}{31352832}}{\tau}^{9}+{\frac {213841
}{17418240}}{\tau}^{7}-{\frac {4069}{7464960}}{\tau}^{5}$\\[1ex]
& $-{\frac {6101
}{7464960}}{\tau}^{3}+{\frac {571}{2488320}}\tau$\\ [1ex]
 8 & ${\frac {1}{7142567040}{\tau}^{26}
}-{\frac {163}{7142567040}}{\tau}^{24}+{\frac {35047}{
23808556800}}{\tau}^{22}-{\frac {41279}{850305600}}{\tau}^{20}+{\frac {
168240799}{190468454400}}{\tau}^{18}$ \\ [1ex]
 & $-{\frac {26988481}{
3023308800}}{\tau}^{16}+{\frac {1266722761}{26453952000}}{\tau}^{14}-{\frac {
3211451701}{26453952000}}{\tau}^{12}+{\frac {21361155917}{193995648000}}{\tau}^{10
}$ \\ [1ex]
& $-{\frac {419199161}{38799129600}}{\tau}^{8}+{\frac 
{889093}{2771366400}}{\tau}^{6}+{\frac {510887}{
923788800}}{\tau}^{4}-{\frac {2016373}{3695155200}}{\tau}^{2}-{\frac{3184811}{3695155200}}$ \\[-1.5ex]
 & \\\hline
\end{tabular}
\end{center}
\caption{The coefficients $C_n(\tau)$ for $0\leq n \leq 8$.}
\label{table1}
\end{table*}


\section{Proof of Proposition \ref{prop1}}\label{section3}

The new asymptotic expansions \eqref{eq7} and \eqref{eq8} were proved under the assumptions that $a$ is large and $\tau$ is bounded. 
Here we show that the restriction on $\tau$ can be relaxed and the expansions \eqref{eq7} and \eqref{eq8} are still valid in the sense of 
generalised asymptotic expansions.

Suppose that $|\tau| = \mathcal{O}\left(|a|^\mu\right)$ as $a\to \infty$ with some constant $\mu < \frac{1}{6}$. With these provisos, 
since $C_n(\tau)$ is a polynomial in $\tau$ of degree $3n+2$, the sequences $C_n (\tau )a^{ - n/2}$ and $\left(-\im\right)^n C_n (\im\tau )a^{ - n/2}$ 
($n=0,1,2,\ldots$) are asymptotic sequences as $a \to \infty$. Furthermore, the remainder terms in \eqref{eq22} and \eqref{eq9} have the same 
order of magnitude as the corresponding first neglected terms $C_N (\tau )a^{ - N/2}$ and $\left(-\im\right)^N C_N (\im\tau )a^{ - N/2}$, respectively. 
This proves that, with the above assumptions on $\tau$, the expansions \eqref{eq7} and \eqref{eq8} are indeed generalised asymptotic expansions.


\section{Proof of Theorems \ref{thminv} and \ref{thmzero}}\label{section4}

We begin with the proof of Theorem \ref{thminv}. The equations \eqref{eq7}, \eqref{inveq} and \eqref{qdef} imply that
\begin{equation}\label{Eeq}
E(\tau):=\sqrt {\frac{\pi }{2}} \exp \left( \frac{{\tau ^2 }}{2}\right)\left( \erfc\left(2^{ - \frac{1}{2}} \tau _0\right) -
\erfc\left(2^{ - \frac{1}{2}} \tau \right) \right) \sim \sum_{n = 0}^\infty \frac{C_n (\tau )}{a^{(n + 1)/2}},
\end{equation}
as $a\to+\infty$. The function $E(\tau)$ satisfies the second order linear homogeneous differential equation
\begin{equation}\label{Ediffeq}
E''(\tau ) - \tau E'(\tau ) - E(\tau ) = 0.
\end{equation}
Expanding $E(\tau)$ into a power series around $\tau=\tau_0$, and substituting it into \eqref{Ediffeq}, we find that
\begin{equation}\label{Eexp}
E(\tau ) = \sum_{k = 1}^\infty  P_k (\tau _0 )\left(\tau  - \tau _0 \right)^k 
\end{equation}
where $P_1 (\tau _0 ) = 1$, $P_2 (\tau _0 ) = \frac{1}{2}\tau _0$ and
\[
k P_{k} (\tau _0 ) = \tau _0 P_{k - 1} (\tau _0 ) + P_{k - 2} (\tau _0 ),
\]
for $k\geq 3$.

We seek for a solution of the asymptotic equality \eqref{Eeq} in the form
\begin{equation}\label{tauexp}
\tau  \sim \tau _0  + \sum_{k = 0}^\infty  \frac{d_k (\tau _0 )}{a^{(k + 1)/2} },\quad a\to +\infty.
\end{equation}
Substituting \eqref{tauexp} into \eqref{Eexp} and expanding in powers of $a^{-1/2}$, we deduce
\[
E(\tau ) \sim \sum_{k = 0}^\infty \left( \sum_{m = 1}^{k + 1} P_m (\tau _0 )\textbf{\textsf{B}}_{k + 1,m} (d_0 (\tau _0 ), \ldots ,d_{k - m + 1} (\tau _0 ))  \right)\frac{1}{a^{(k + 1)/2} } ,
\]
as $a\to+\infty$, where $\textbf{\textsf{B}}_{k,m}$ are the partial ordinary Bell polynomials (see, e.g., \cite[Appendix]{Nemes13}). Similarly, employing \eqref{eq10a}, we obtain
\begin{align*}
\sum_{n = 0}^\infty \frac{C_n (\tau )}{a^{(n + 1)/2} } 
& \sim \sum_{n = 0}^\infty  \sum_{m = 0}^{3n + 2} \frac{c_{n,m}}{a^{(n - m + 1)/2} }\left( \frac{\tau }{a^{1/2} } \right)^m  \\ 
& \sim \sum_{k = 0}^\infty  \left( \sum_{n = 0}^k \sum_{m = 0}^{3n + 2} c_{n,m} 
\textbf{\textsf{B}}_{k - n+m,m} \left(\tau _0 ,d_0 (\tau _0 ), \ldots ,d_{k - n - 1} (\tau _0 )\right)  \right)\frac{1}{a^{(k + 1)/2}} ,
\end{align*}
as $a\to+\infty$. By equating to zero the coefficients of powers of $a^{-1/2}$, and using $P_1 (\tau _0 ) = 1$, 
we derive the recurrence relation $d_0 (\tau _0 ) = \frac{1}{3}\tau _0^2  - \frac{1}{3}$ and
\begin{gather}\label{drec}
\begin{split}
d_k (\tau _0 ) = & - \sum_{m = 2}^{k + 1} P_m (\tau _0 )\textbf{\textsf{B}}_{k + 1,m} (d_0 (\tau _0 ), \ldots ,d_{k - m + 1} (\tau _0 )) \\ & + \sum_{n = 0}^k \sum_{m = 0}^{3n + 2} c_{n,m} \textbf{\textsf{B}}_{k - n+m,m} (\tau _0 ,d_0 (\tau _0 ), \ldots ,d_{k - n - 1} (\tau _0 )) ,
\end{split}
\end{gather}
for $k\geq 1$. We can use this recurrence and induction on $k$ to prove that $d_k ( - \tau _0 ) = \left( - 1\right)^k d_k (\tau _0 )$, i.e., 
that the even and odd-order polynomials are even and odd functions, respectively. The identity clearly holds for $k=0$. 
Suppose that $k\geq 1$. Using induction, it is easy to see that $P_m ( - \tau _0 ) = \left( - 1\right)^{m + 1} P_m (\tau _0 )$. 
From the properties of the partial ordinary Bell polynomials, we can infer that
\begin{align*}
\textbf{\textsf{B}}_{k + 1,m} (d_0 ( - \tau _0 ), \ldots ,d_{k - m + 1} ( - \tau _0 )) & = \textbf{\textsf{B}}_{k + 1,m} (d_0 (\tau _0 ), \ldots ,\left(-1\right)^{k - m + 1} d_{k - m + 1} (\tau _0 )) \\ & = \left(-1\right)^{k - m + 1} \textbf{\textsf{B}}_{k + 1,m} (d_0 (\tau _0 ), \ldots ,d_{k - m + 1} (\tau _0 ))
\end{align*}
and
\begin{align*}
& \textbf{\textsf{B}}_{k - n + m,m} ( - \tau _0 ,d_0 ( - \tau _0 ), \ldots ,d_{k - n - 1} ( - \tau _0 )) \\ & = \textbf{\textsf{B}}_{k - n + m,m} ( - \tau _0 ,d_0 (\tau _0 ), \ldots ,\left(-1\right)^{k - n - 1} d_{k - n - 1} (\tau _0 ))\\ & = \left(-1\right)^{k - n + m} \textbf{\textsf{B}}_{k - n + m,m} (\tau _0 ,d_0 (\tau _0 ), \ldots ,d_{k - n - 1} (\tau _0 )) .
\end{align*}
Finally, $C_n ( - \tau ) = \left( - 1\right)^n C_n (\tau )$ implies $c_{n,m}  = \left( - 1\right)^{ - n + m} c_{n,m}$. 
Employing all these relations in \eqref{drec}, it follows readily that $d_k ( - \tau _0 ) = \left( - 1\right)^k d_k (\tau _0 )$.

A similar argument would show that the degree of $d_k(\tau_0)$ as a polynomial in $\tau_0$ is at most $3k+2$. To prove that the degree is actually $k+2$, we need another recurrence for these polynomials. In order to obtain this recurrence, we proceed as follows. If we differentiate both sides of the asymptotic equality
\[
\tfrac{1}{2}\erfc\left(2^{ - \frac{1}{2}} \tau _0 \right) \sim Q\left( a,a + \tau _0 a^{\frac{1}{2}}  + \sum_{k = 0}^\infty  \frac{d_k (\tau _0 )}{a^{k/2}} \right),\quad a\to +\infty,
\]
with respect to $\tau_0$, we find
\begin{align*}
 - \frac{1}{\sqrt {2\pi }}\exp \left(  - \frac{\tau_0 ^2 }{2} \right) \sim - & \frac{1}{\sqrt {2\pi }}\frac{1}{\Gamma ^ *  (a)}\left( 1 + \left( \tau _0  + \sum_{k = 0}^\infty  \frac{d_k (\tau _0 )}{a^{(k + 1)/2}} \right)a^{ - \frac{1}{2}}  \right)^{a - 1} \\ & \times \exp \left(  - \tau _0 a^{\frac{1}{2}}  - \sum_{k = 0}^\infty  \frac{d_k (\tau _0 )}{a^{k/2} }  \right)\left( 1 + \sum_{k = 0}^\infty  \frac{d'_k (\tau _0 )}{a^{(k + 1)/2}} \right)
\end{align*}
as $a\to+\infty$, where $\Gamma ^ *  (a): = \Gamma (a)\e^a a^{1/2 - a} (2\pi )^{ - 1/2}$ is the so-called scaled gamma function. Taking logarithms, we can assert that
\begin{align*}
 - \frac{\tau_0 ^2}{2} \sim & - \log \Gamma ^ *  (a) + (a - 1)\log \left( 1 + \left( \tau _0  + \sum_{k = 0}^\infty  \frac{d_k (\tau _0 )}{a^{(k + 1)/2} } \right)a^{ - \frac{1}{2}} \right) \\ & - \tau _0 a^{\frac{1}{2}}  - \sum_{k = 0}^\infty  \frac{d_k (\tau _0 )}{a^{k/2} }  + \log \left( 1 + \sum_{k = 0}^\infty \frac{d'_k (\tau _0 )}{a^{(k + 1)/2} } \right).
\end{align*}
Expanding the right-hand side in powers of $a^{-1/2}$, applying the known asymptotic expansion of $\log \Gamma^*(a)$ 
\cite[\href{http://dlmf.nist.gov/5.11.1}{Eq.~5.11.1}]{NIST:DLMF}, and comparing the coefficients of powers of $a^{-1/2}$ on each side, we deduce
\begin{gather}\label{altdrec}
\begin{split}
& \sum_{m = 2}^{k + 2} \frac{\left(-1\right)^m }{m}\textbf{\textsf{B}}_{k + 2,m} (\tau _0 ,d_0 (\tau _0 ), \ldots ,d_{k - m + 1} (\tau _0 )) \\ & + \sum_{m = 1}^k \frac{{\left(-1\right)^m }}{m}\left( \textbf{\textsf{B}}_{k,m} (d'_0 (\tau _0 ), \ldots ,d'_{k - m} (\tau _0 )) - \textbf{\textsf{B}}_{k,m} (\tau _0 ,d_0 (\tau _0 ), \ldots ,d_{k - m - 1} (\tau _0 )) \right) \\ & = \begin{cases} -\cfrac{{B_{2n} }}{{2n(2n - 1)}} &\mbox{if } k=4n-2, \\ 
0 & \mbox{otherwise.}\end{cases}
\end{split}
\end{gather}
Here $B_{2n}$ denotes the even-order Bernoulli numbers \cite[\href{http://dlmf.nist.gov/24.2.i}{\S24.2(i)}]{NIST:DLMF}. We can extract out the highest coefficient $d_{k - 1} (\tau _0 )$ and write the above expression as
\begin{align*}
& \tau _0 d_{k - 1} (\tau _0 ) - d'_{k - 1} (\tau _0 ) + \tfrac{1}{2}\sum_{m = 1}^{k - 1} {d_{m - 1} (\tau _0 )d_{k - m - 1} (\tau _0 )} \\ & + \sum_{m = 3}^{k + 2} {\frac{{\left(-1\right)^m }}{m}\textbf{\textsf{B}}_{k + 2,m} (\tau _0 ,d_0 (\tau _0 ), \ldots ,d_{k - m + 1} (\tau _0 ))} 
\\ & + \sum_{m = 2}^k {\frac{{\left(-1\right)^m }}{m}\textbf{\textsf{B}}_{k,m} (d'_0 (\tau _0 ), \ldots ,d'_{k - m} (\tau _0 ))}  - \sum_{m = 1}^k {\frac{{\left(-1\right)^m }}{m}\textbf{\textsf{B}}_{k,m} (\tau _0 ,d_0 (\tau _0 ), \ldots ,d_{k - m - 1} (\tau _0 ))} 
\\ & = \begin{cases} -\cfrac{{B_{2n} }}{{2n(2n - 1)}} &\mbox{if } k=4n-2, \\ 
0 & \mbox{otherwise.}\end{cases}
\end{align*}
Now a simple induction argument using the properties of the partial ordinary Bell polynomials shows that the degree of $d_{k-1}(\tau_0)$ as a polynomial in $\tau_0$ is at most $k+1$. It remains to prove that the leading coefficients of these polynomials are non-zero. Let us denote $\delta _1  = 1$ and
\[
d_k (\tau _0 ) = \delta _{k + 2} \tau _0^{k + 2}  +  \cdots ,
\]
for $k\geq 0$. From \eqref{altdrec}, we can infer that
\[
\sum_{m = 2}^k \frac{\left(-1\right)^m}{m}\textbf{\textsf{B}}_{k,m} (\delta _1 , \ldots ,\delta _{k - m + 1} )  = 0,\quad k \ge 3,
\]
or, equivalently,
\[
\sum_{m = 1}^k \frac{\left(-1\right)^{m + 1}}{m}\textbf{\textsf{B}}_{k,m} (\delta _1 , \ldots ,\delta _{k - m + 1} )  = \delta _k ,\quad k \ge 3.
\]
By the definition of the partial ordinary Bell polynomials, this equality translates to the following relation between formal power series
\[
\log \left( 1 + \sum_{k = 1}^\infty \delta _k x^k   \right) = -\frac{x^2}{2} + \sum_{k = 1}^\infty \delta _k x^k .
\]
By the results of the paper \cite{Brassesco2011}, we have the explicit expression
\[
\delta _k  = \frac{1}{k!}\left[ \frac{\d^{k - 1} }{\d x^{k - 1} }\left( \frac{x^2/2}{\e^x  - x - 1} \right)^{k/2}  \right]_{x = 0} 
\]
provided $k \ge 3$. Combining this expression with the results of the paper \cite{Nemes16}, we find that the sequence $\delta_k$ is related to the coefficients of the large-$a$ asymptotic expansion of $\Gamma (a,a)$ as follows:
\[
\Gamma (a,a) \sim \sqrt {\frac{\pi }{2}} a^{a - 1/2} \e^{ - a} \left( 1 - \frac{1}{3}\sqrt {\frac{2}{\pi a}}  + \sum_{k = 2}^\infty  \frac{(k + 1)!\delta _{k + 1} }{2^{k/2} \Gamma \left( \frac{k}{2} + 1 \right)}\frac{1}{a^{k/2}} \right), \quad a\to +\infty.
\]
The coefficients in this asymptotic expansion are known to be non-zero as a consequence of their integral representations given in 
\cite{Nemes16}. Thus, $\delta _k\neq 0$ if $k\geq 3$. The case $k=2$ is trivial since $\delta_2=\frac{1}{3}$.

We continue with the proof of Theorem \ref{thmzero}. A simple manipulation of \eqref{eq8} shows that
\begin{equation}\label{eq15}
\left(-x\right)^{a} \gamma^\ast (a,x) \sim \cos (\pi a) + \frac{\sin (\pi a)}{\pi}\left( 2\sqrt \pi  F\left(2^{ - \frac{1}{2}} \tau\right) 
+ \sqrt {\frac{2\pi}{-a}} \sum_{n = 0}^\infty \frac{\left(-\im\right)^n C_n (\im\tau )}{\left(-a\right)^{n/2}}  \right)\exp \left( \frac{\tau ^2 }{2} \right),
\end{equation}
as $a\to -\infty$, uniformly with respect to bounded real values of $\tau$ (compare \cite[\href{http://dlmf.nist.gov/8.11.11}{Eq.~8.11.11}]{NIST:DLMF}). 
Here $x=a-\tau \left(-a\right)^{\frac{1}{2}}$. Now, let $x_{-}(a)$ be the unique negative zero of $\gamma^\ast (a,x)$ and write 
$x_{-}(a)=a-\tau \left(-a\right)^{\frac{1}{2}}$. Employing the definition \eqref{eq12} of $\tau_1$, \eqref{eq15} may be rewritten in the form
\[
\sqrt 2\exp \left(  - \frac{\tau^2}{2} \right)\left( F\left(2^{ - \frac{1}{2}} \tau _1\right)\exp \left( \frac{\tau _1^2}{2} \right) - F\left(2^{ - \frac{1}{2}} \tau\right)
\exp \left( \frac{\tau ^2}{2} \right) \right) \sim \sum_{n = 0}^\infty \frac{\left(-\im\right)^n C_n (\im\tau )}{\left( - a\right)^{(n + 1)/2}} ,
\]
as $a\to -\infty$. This asymptotic equality is further equivalent to
\begin{equation}\label{invaseq}
\sqrt {\frac{\pi }{2}} \exp \left(  - \frac{\tau ^2}{2} \right)\left( \erfc(2^{ - \frac{1}{2}} \im\tau _1 ) - \erfc(2^{ - \frac{1}{2}} \im\tau ) \right) \sim \sum_{n = 0}^\infty \frac{( - \im)^{n + 1} C_n (\im\tau )}{( - a)^{(n + 1)/2}},
\end{equation}
as $a\to -\infty$. Thus, the inversion problem is formally equivalent to that of \eqref{Eeq}, with $\im\tau$, $\im\tau_1$ and $ - a\e^{\pi\im}$ 
in place of $\tau$, $\tau_0$ and $a$. Therefore, if we seek for a solution of the asymptotic equality \eqref{invaseq} in the form
\[
\tau  \sim \tau _1  + \sum\limits_{k = 0}^\infty  \frac{\widetilde d_k (\tau _1 )}{( - a)^{(k + 1)/2} },\quad a\to -\infty,
\]
using the above procedure, we find that $\widetilde d_k (\tau _1 ) =  - ( - \im)^k d_k (\im\tau _1 )$.

We close this section by commenting on the implementation of the recurrence relation \eqref{drec}. To evaluate the partial ordinary Bell polynomials, we use the formal generating function \cite[Appendix]{Nemes13}
\begin{equation}\label{eq19a}
\exp \left( y\sum\limits_{n = 1}^\infty \alpha _n x^n  \right) = \sum\limits_{k = 0}^\infty \sum\limits_{m = 0}^\infty  \frac{\textbf{\textsf{B}}_{k,m} \left( \alpha _1 ,\alpha _2 , \ldots ,\alpha _{k - m + 1} \right)}{m!}x^k y^m.
\end{equation}
The Taylor coefficients of this generating function can be easily obtained with any computer algebra package. We combine \eqref{eq19a} with \eqref{drec} to compute the $d_{k}(\tau_0)$'s. The explicit forms of the first several coefficients are presented in Table \ref{table2} below. Note that when evaluating $\textbf{\textsf{B}}_{k,m}$ ($m=0,1,\ldots,k$) via \eqref{eq19a}, the sum $\displaystyle \sum_{n=1}^\infty$ can be replaced by the finite sum $\displaystyle \sum_{n=1}^{k+1}$.

\begin{table*}[!ht]
\begin{center}
\begin{tabular}
[c]{l l}  \hline
& \\ [-2ex]
 $n$ & $d_n(\tau_0)$ \\ [-2ex]
& \\ \hline 
 & \\ [-1.5ex]
 0 & ${\frac{1}{3}}\tau_{0}^{2}-{\frac{1}{3}}$\\ [1ex]
 1 & ${\frac{1}{36}}\tau_{0}^{3}-{\frac{7}{36}}\tau_{0}$\\ [1ex]
 2 & $-{\frac{1}{270}}\tau_{0}^{4}-{\frac{7}{810}}\tau_{0}^{2}+{\frac{8}{405}}$\\[1ex]
 3 & ${\frac{1}{4320}}\tau_{0}^{5}+{\frac{8}{1215}}\tau_{0}^{3}-{\frac{433}{38880}}\tau_{0}$\\[1ex]
 4 & ${\frac{1}{17010}}\tau_{0}^{6}-{\frac{1}{840}}\tau_{0}^{4}-{\frac{923}{204120}}\tau_{0}^{2}+{\frac{184}{25515}}$\\[1ex]
 5 & $-{\frac{139}{5443200}}\tau_{0}^{7}-{\frac{1451}{48988800}}\tau_{0}^{5}+{\frac{289517}{146966400}}\tau_{0}^{3}+{\frac{289717}{146966400}}\tau_{0}$ \\[1ex]
 6 & ${\frac{1}{204120}}\tau_{0}^{8}+{\frac{769}{9185400}}\tau_{0}^{6}-{\frac{151}{874800}}\tau_{0}^{4}-{\frac{104989}{55112400}}\tau_{0}^{2}+{\frac{2248}{3444525}}$ \\[1ex]
 7 & $-{\frac{571}{2351462400}}\tau_{0}^{9}-{\frac{1087}{41990400}}\tau_{0}^{7}-{\frac{30469}{235146240}}\tau_{0}^{5}+{\frac{219257}{661348800}}\tau_{0}^{3}+{\frac{1500053}{846526464}}\tau_{0}$\\ [1ex]
 8 & $-{\frac{281}{1515591000}}\tau_{0}^{10}+{\frac{49271}{15588936000}}\tau_{0}^{8}+{\frac{997903}{15588936000}}\tau_{0}^{6}+{\frac{101251277}{654735312000}}\tau_{0}^{4}$ \\ [1ex]
 & $-{\frac{96026707}{280600848000}}\tau_{0}^{2}-{\frac{19006408}{15345358875}}$ \\[1ex]
 9 & ${\frac{163879}{2172751257600}}\tau_{0}^{11}+{\frac{209488529}{293321419776000}}\tau_{0}^{9}-{\frac{252836779}{20951529984000}}\tau_{0}^{7}-{\frac{15974596457}{146660709888000}}\tau_{0}^{5}$ \\[1ex]
 & $-{\frac{556030221167}{2639892777984000}}\tau_{0}^{3}+{\frac{487855454729}{2639892777984000}}\tau_{0}$ \\[1ex]
 10 & $-{\frac{5221}{354648294000}}\tau_{0}^{12}-{\frac{1316963}{2708223336000}}\tau_{0}^{10}-{\frac{329677}{255346771680}}\tau_{0}^{8}+{\frac{1386098437}{53622822052800}}\tau_{0}^{6}$ \\ [1ex]
 & $+{\frac{11414859619}{71497096070400}}\tau_{0}^{4}+{\frac{1069525622411}{3217369323168000}}\tau_{0}^{2}-{\frac{5667959576}{12567848918625}}$ \\[-1.5ex]
 & \\\hline
\end{tabular}
\end{center}
\caption{The coefficients $d_n(\tau_0)$ for $0\leq n \leq 10$.}
\label{table2}
\end{table*}


\section{Numerical comparisons}\label{section5}
In this section, we compare the global uniform asymptotic expansion \eqref{eq5} with our new transition region expansion \eqref{eq7}.
We denote $\lambda=z/a$ as usual.

\subsection{The transition region expansion}
According to Proposition \ref{prop1}, the transition region expansion \eqref{eq7} holds for $\tau=\mathcal{O}\left(|a|^\mu\right)$, as $a\to\infty$,
with $\mu<\frac16$, that is for $\lambda-1=\mathcal{O}\left(|a|^\nu\right)$, with $\nu<-\frac13$. The main reason for this is that the coefficients
$C_n(\tau)$ are polynomials of degree $3n+2$, and if $\tau\sim K a^{1/6}$ held, all the terms in the divergent asymptotic expansion \eqref{eq7} would be of the same order. It follows from \eqref{eq10a} and \eqref{eq10c} that if we replaced the $C_n(\tau)$'s in \eqref{eq7} by their leading terms, then the resulting infinite series would actually converge. Numerical experiments illustrate that even in the case that $\tau\sim K a^{1/6}$, the expansion \eqref{eq7} still provides reasonable approximations. Accordingly, the expansion \eqref{eq7} has an asymptotic property when $\tau=\mathcal{O}\left(|a|^\mu\right)$, with $\mu<\frac16$, and we can allow $\mu$ to be close to $\frac16$.

\subsection{The global uniform asymptotic expansion}
The global uniform asymptotic expansion \eqref{eq5} has the advantage that it is valid as $a\to\infty$, uniformly with respect to $\lambda$ in the sector $|\arg \lambda|\leq 2\pi-\delta< 2\pi$. The main drawback of the expansion is that the coefficients are significantly more complicated than the $b_n(\lambda)$'s and $C_n(\tau)$'s above, which are polynomials satisfying simple recurrence relations. The evaluation of the coefficients $c_n(\eta)$ near the point $\eta=0$ ($\lambda=1$) is especially difficult because they possess a removable singularity at this point. In Appendix \ref{appendixb}, we give several methods to compute the $c_n(\eta)$'s in a stable manner.

\begin{figure}[htbp]
\centering
\leavevmode
\hbox{\placefigure{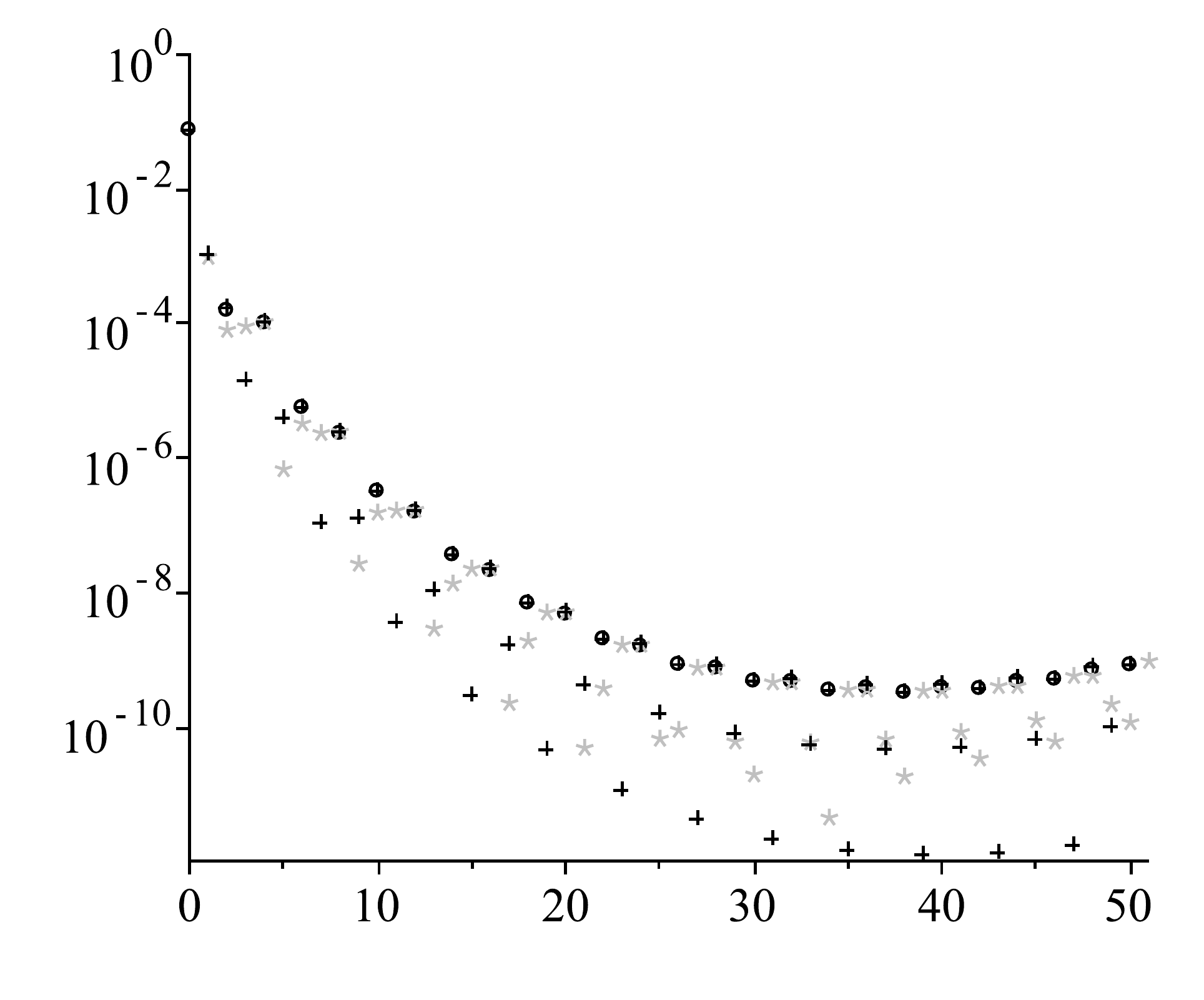}{0.45\hsize}%
\placetext{$a=3$}{0.85}{0.5}%
\placetext{$\tau=0.1$}{0.85}{0.44}%
\placetext{$n$}{1.02}{0.053}%
\hbox to \wd0{\hss}}
\caption{The absolute value of the $n$th term in \eqref{eq7} (crosses), the absolute value of the $\frac{1}{2}n$th term (for even $n$) in \eqref{eq5} (circles), and the absolute value of the remainder in \eqref{eq7} after $n$ terms (grey asterisks).}
\label{fig1}
\end{figure}

In Figure \ref{fig1}, we compare the terms of the asymptotic series \eqref{eq7} and \eqref{eq5} in the case that $a=3$ and $\tau=0.1$.
Since $\tau$ it small, it is illustrated that $C_{2n}(\tau)\approx c_{n}(\eta)$, which follows from the identity $C_{2n}(0)=c_{n}(0)$.
The terms first decrease in magnitude, reach a minimum, and then increase. This is the typical behaviour of Gevrey-1 divergent series (see, e.g., \cite{MS16}). The grey asterisks in Figure \ref{fig1} are the absolute values of the remainders in \eqref{eq7} after $n$ terms. As expected, they are approximately of the same size as the corresponding first neglected terms. This will also be the case for the uniform asymptotic expansion \eqref{eq5}. 
It is suggested by Figure \ref{fig1} that in the case of $n=34$, our approximation produces 11 correct digits whereas the uniform asymptotic approximation \eqref{eq5}
gives us only 9 correct digits.

\begin{figure}[htbp]
\centering
\leavevmode
\hbox{\placefigure{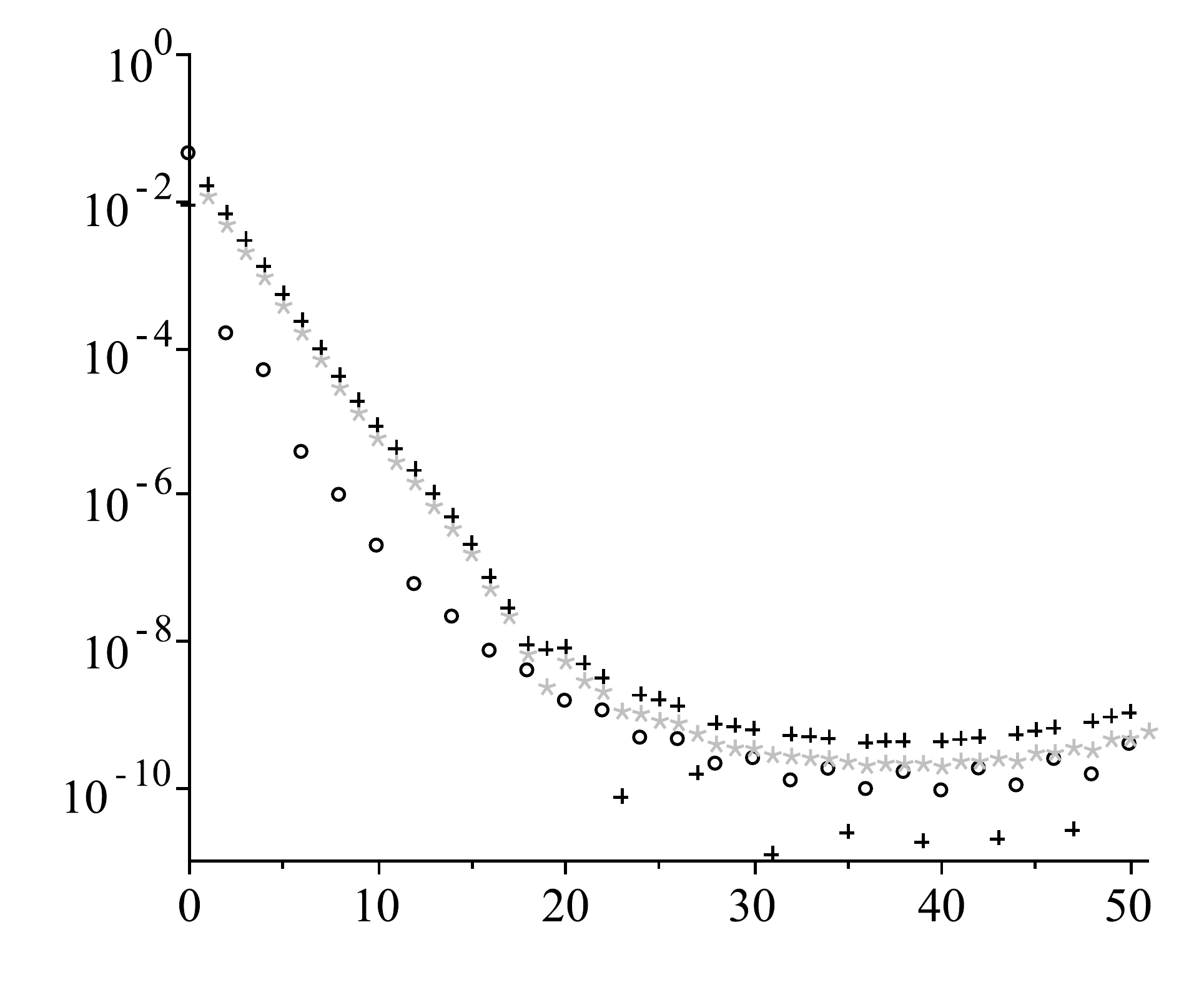}{0.45\hsize}%
\placetext{$a=3$}{0.85}{0.5}%
\placetext{$\tau=1.1$}{0.85}{0.44}%
\placetext{$n$}{1.02}{0.05}%
\hbox to \wd0{\hss}}
\caption{The absolute value of the $n$th term in \eqref{eq7} (crosses), the absolute value of the $\frac{1}{2}n$th term (for even $n$) in \eqref{eq5} (circles), and the absolute value of the remainder in \eqref{eq7} after $n$ terms (grey asterisks).}
\label{fig2}
\end{figure}

Similarly, in Figure \ref{fig2}, we consider the case that $a=3$ and $\tau=1.1$. No significant difference is seen between Figures \ref{fig1} and \ref{fig2}.

\begin{figure}[htbp]
\centering
\leavevmode
\hbox{\placefigure{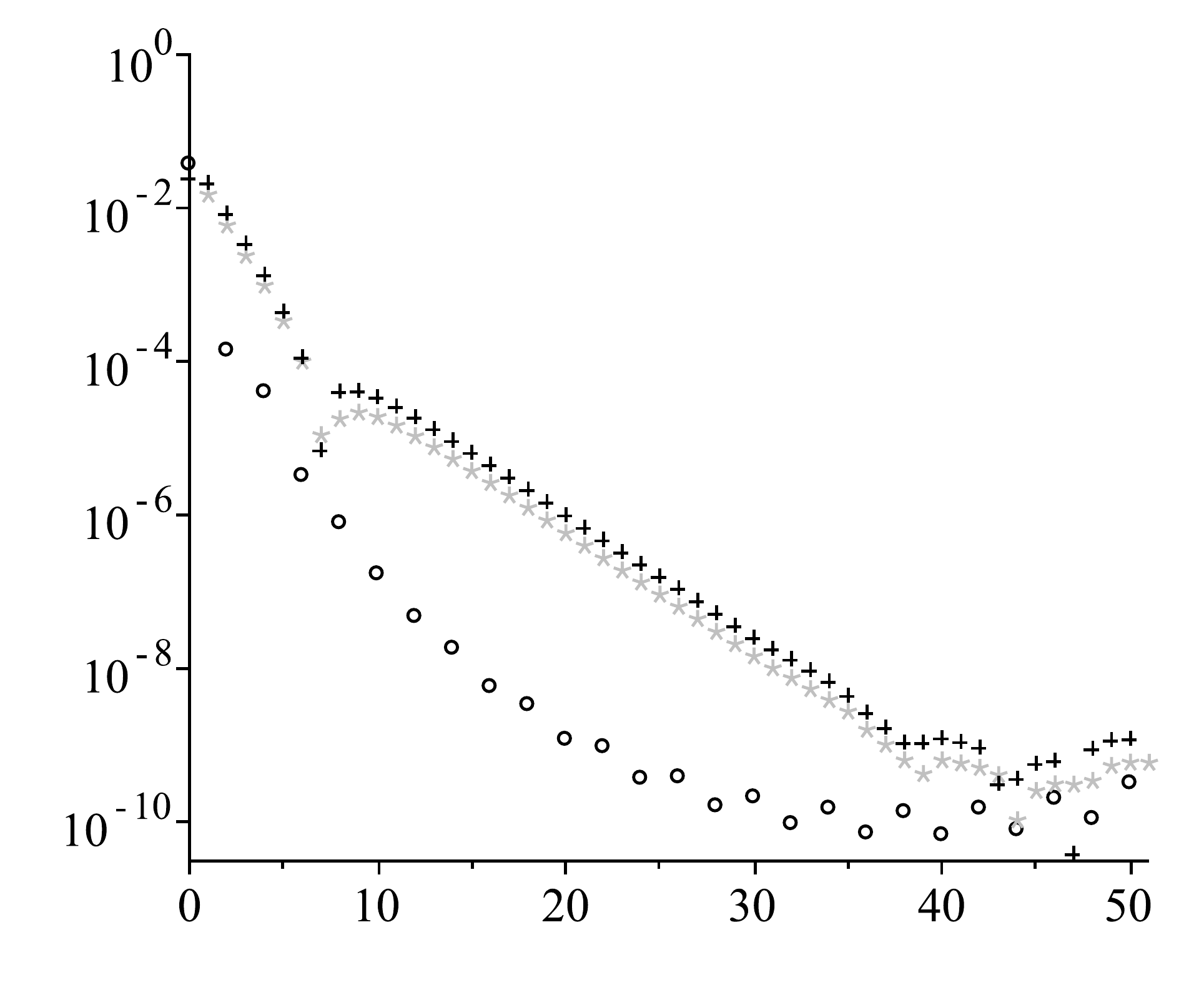}{0.45\hsize}%
\placetext{$a=3$}{0.8}{0.5}%
\placetext{$\tau=1.321$}{0.8}{0.44}%
\placetext{$n$}{1.02}{0.047}%
\hbox to \wd0{\hss}}
\caption{The absolute value of the $n$th term in \eqref{eq7} (crosses), the absolute value of the $\frac{1}{2}n$th term (for even $n$) in \eqref{eq5} (circles), and the absolute value of the remainder in \eqref{eq7} after $n$ terms (grey asterisks).}
\label{fig3}
\end{figure}

Finally, in Figure \ref{fig3}, we consider the case when $a=3$ and $\tau=1.321$, that is, when $\tau\approx 1.1 a^{1/6}$. Thus, the divergent series \eqref{eq7} does not have an asymptotic property anymore. Although the behaviour of the terms in \eqref{eq7} has clearly changed, it seems that the magnitudes of the smallest terms in \eqref{eq7} and \eqref{eq5} remained approximately equal.

\section*{Acknowledgement} The authors thank the anonymous referee for his/her helpful comments and suggestions on the manuscript.

\appendix

\section{Asymptotic expansions away from the transition points}\label{appendixa}

The outer expansions \eqref{eq1} and \eqref{eq2} break down near the transition point, that is, when 
$\lambda$ approaches $1$, since their terms become singular in this limit. Note that in Theorem \ref{thm1} we take $z=a+\tau a^{1/2}$, thus $\lambda=z/a=1+\tau a^{-1/2}$. In the forthcoming we examine the $\tau$-region of validity of these outer expansions.

For \eqref{eq1} and \eqref{eq2} to have an asymptotic property we need $\left(z-a\right)^2/a=a\left(\lambda-1\right)^2$ to be large as $a\to\infty$. 
Hence, we require that $\lambda-1=\mathcal{O}\left(|a|^\nu\right)$, as $a\to\infty$, with $\nu>-\frac12$. Therefore, since $\tau=(\lambda-1)a^{1/2}$, we require that $\tau=\mathcal{O}\left(|a|^\mu\right)$, as $a\to\infty$, with $\mu>0$.

We can also consider the optimal number of terms of the divergent series \eqref{eq1} and \eqref{eq2} to illustrate the necessity of the condition $\nu>-\frac12$. For this, we need the asymptotic behaviour of the coefficients $b_n(\lambda)$ as $n\to+\infty$. We take $\lambda>0$, $\lambda \neq 1$. Then, according to \cite[Eq. 23, p. 161]{Dingle73}, \cite{Nemes16}, we have
\begin{equation}\label{lateterm1}
b_n (\lambda)=\frac{1}{\sqrt {2\pi}}\frac{\Gamma \left( n+\frac{1}{2} \right)\left(\lambda-1\right)^{2n+1} }{\left(\lambda-\log\lambda-1\right)^{n + \frac{1}{2}} } \operatorname{sgn} (\lambda  - 1) \left( 1 +
\mathcal{O}_\lambda\left( \frac{1}{n} \right)\right)
\end{equation}
as $n \to +\infty$. Now let $n$ be the optimal number of terms of the divergent series \eqref{eq1} and \eqref{eq2}, i.e., the index of the numerically least term. 
Then, by \eqref{lateterm1}, approximately
\[
1\approx\left|\frac{\frac{\left(-a\right)^{n+1}b_{n+1}(\lambda)}{\left(z-a\right)^{2n+3}}}{\frac{\left(-a\right)^{n}b_{n}(\lambda)}{\left(z-a\right)^{2n+1}}}\right|
=\left|\frac{b_{n+1}(\lambda)}{a\left(\lambda-1\right)^2b_{n}(\lambda)}\right|\sim\left|\frac{n}{a(\lambda-1-\log\lambda)}\right|,
\]
as $n\to+\infty$. In the case that $\lambda-1=\mathcal{O}\left(|a|^\nu\right)$, with $-\frac12<\nu<0$, we have $2n\approx |a|^{1+2\nu}$, and the optimal
number of terms shrinks to zero as $\nu$ approaches $-\frac12$.

\section{The coefficients in the uniform asymptotic expansions}\label{appendixb}

As mentioned above, in Temme \cite{Temme1979} power series expansions in $\eta$ are given for the coefficients $c_n(\eta)$. 
Here we consider the same expansion, but also one in powers of $\lambda-1$; the latter seems more natural, and the details are simpler. In both cases, we give (new) recurrence relations for the coefficients. We introduce the expansions
\begin{equation}\label{eq6b}
c_n(\eta)=\sum_{k=0}^\infty e_{k,n}\left(\lambda-1\right)^k=\sum_{k=0}^\infty f_{k,n}\eta^k.
\end{equation}
Observe that $c_0(\eta)$ satisfies the equation
\[
2\left(\lambda-\log\lambda-1\right)\left((\lambda-1)c_0(\eta)-1\right)^2=\left(\lambda-1\right)^2,
\]
and since $\lambda\eta\frac{\d\eta}{\d\lambda}=\lambda -1$, the recurrence relation in \eqref{eq6a} can be written as
\[
c_n (\eta) = \frac{\gamma _n }{\lambda  - 1} + \frac{\lambda}{\lambda-1}\frac{\d c_{n - 1} (\eta)}{\d\lambda},\quad n \ge 1.
\]
Therefore, we obtain for the coefficient $e_{k,n}$ the recurrence relations
\begin{equation}\label{eq6e}
e_{k,0}=\frac{\left(-1\right)^{k+1}}{k+3}-2\sum_{l=1}^{k}\frac{\left(-1\right)^{l}}{l+2}e_{k-l,0}
-\sum_{l=1}^{k}\sum_{m=1}^{l} \frac{\left(-1\right)^{m}}{m+1}e_{k-l,0}e_{l-m,0},
\end{equation}
thus $e_{0,0}=-\frac13$, and
\begin{equation}\label{eq6f}
e_{k,n+1}=(k+1)e_{k+1,n}+(k+2)e_{k+2,n}.
\end{equation}
Hence, it is straightforward to compute the first several coefficients. Note that if we know the first $K$ Taylor coefficients of $c_n(\eta)$, then we can use \eqref{eq6f} to compute the first $K-2$ Taylor coefficients of $c_{n+1}(\eta)$. In each step we lose two Taylor coefficients. Fortunately, this is not a serious problem, since it is easy to compute many of the coefficients for $c_0(\eta)$ via \eqref{eq6e}.

Regarding the coefficients $f_{k,n}$ in \eqref{eq6b}, we observe that  $c_0(\eta)$ satisfies the differential equation
\[
\eta \frac{\d c_{0} (\eta)}{\d\eta}+\eta^2 c_{0}^3 (\eta)+\left(\eta^2+3\eta\right)c_{0}^2 (\eta)+(2\eta+3)c_{0} (\eta)+1=0,
\]
and that we can rewrite \eqref{eq6a} as
\begin{equation}\label{eq6h}
c_n (\eta) = \gamma _n c_0 (\eta ) + \frac{1}{\eta}\left(\frac{\d c_{n - 1} (\eta )}{\d\eta}+\gamma_n\right).
\end{equation}
Therefore, we obtain for the coefficient $f_{k,n}$ the recurrence relations
\begin{equation}\label{eq6j}
-(k+3)f_{k,0}=2f_{k-1,0}+3\sum_{l=1}^{k}f_{l-1,0}f_{k-l,0}+\sum_{l=2}^{k}f_{l-2,0}f_{k-l,0}+\sum_{l=2}^{k}\sum_{m=0}^{k-l}f_{l-2,0}f_{m,0}f_{k-l-m,0},
\end{equation}
with $f_{0,0}=-\frac13$, and
\begin{equation}\label{eq6i}
f_{k,n+1}=(k+2)f_{k+2,n}-f_{1,n}f_{k,0}.
\end{equation}
Note that from \eqref{eq6h} it follows that $\gamma_{n+1}+f_{1,n}=0$, and we have used this to express the final term in \eqref{eq6i} 
in terms of $f_{1,n}$.

In \eqref{eq6b}, the first series converges for $|\lambda-1|<1$ and the second for $|\eta|<2\sqrt\pi$. 
Hence, although the $\lambda$-sum has the advantage that
it is in terms of the original parameter, the $\eta$-sum has the advantage that it converges in a larger region.

Another way to compute the $c_n(\eta)$'s in a stable manner is to use the integral representation \cite{Dunster1998}
\[
c_n(\eta)=\frac{\im\left(-1\right)^n\Gamma(n+\frac12)}{\left(2\pi\right)^{\frac32}}
\oint_{\{1,\lambda\}}\frac{\,\d t}{(t-\lambda)\left(t-1-\log t\right)^{n+\frac12}},
\]
where the contour of integration is a simple loop surrounding the points $1$ and $\lambda$ in the positive sense. 
Taking $\lambda$ close to $1$, we can choose for the contour a circle with centre $1$ and radius $r<1$. According to the paper \cite{TW14}, 
the trapezoidal rule converges exponentially fast for this type of integral and we obtain a simple method to compute the coefficients via
\begin{equation}\label{eq19d}
c_n(\eta)\approx \frac{\Gamma(n+\frac12)}{2M\sqrt{2\pi}}
\sum_{m=1-M}^M\frac{\sqrt{\displaystyle\frac{\omega_m^2}{\omega_m-\log(\omega_m+1)}}}{(\lambda-\omega_m-1)
\left(\log(\omega_m+1)-\omega_m\right)^n},~~~\textrm{where}~~\omega_m=r\e^{\pi\im m/M}.
\end{equation}
The introduction of the factor $\omega_m^2$ inside the square root makes it single-valued, and therefore computer algebra packages can work with \eqref{eq19d} without any problems.


\begin{thebibliography}{10}

\bibitem{Brassesco2011}
S. Brassesco, M. A. M\'endez, The asymptotic expansion for $n!$ and the Lagrange inversion formula, \emph{Ramanujan J.} \textbf{24} (2011), no. 2, pp. 219--234.

\bibitem{Dingle73}
R. B. Dingle, \emph{Asymptotic Expansions: Their Derivation and Interpretation}, Academic Press, London, 1973.

\bibitem{Dunster1996}
T. M. Dunster, Asymptotics of the generalized exponential integral, and error bounds in the uniform asymptotic smoothing of its Stokes discontinuities, \emph{Proc. Roy. Soc. London Ser. A} \textbf{452} (1996), no. 1949, pp. 1351--1367.

\bibitem{Dunster1998}
T. M. Dunster, R. B. Paris, S. Cang, On the high-order coefficients in the uniform asymptotic expansion for the incomplete gamma function, \emph{Methods Appl. Anal.} \textbf{5} (1998), no. 3, pp. 223--247.

\bibitem{Gautschi1959}
W. Gautschi, Exponential integral $\int_1^\infty \e^{-xt} t^{-n} \d t$ for large values of $n$, \emph{J. Res. Nat. Bur. Standards} \textbf{62} (1959), pp. 123--125.

\bibitem{Giles2016}
M. B. Giles, Algorithm 955: Approximation of the inverse Poisson cumulative distribution function, \emph{ACM Trans. Math. Softw.} \textbf{42} (2016), no. 1, Article 7, 22 pp.

\bibitem{FLP05}
C. Ferreira, J. L. L\'opez, E. P\'erez-Sinus\'{\i}a, Incomplete gamma functions for large values of their variables, \emph{Adv. in Appl. Math.} \textbf{34} (2005), no. 3, pp. 467--485.

\bibitem{Fettis1979}
H. E. Fettis, An asymptotic expansion for the upper percentage points of the $\chi^2$-distribution, \emph{Math. Comp.} \textbf{33} (1979), no. 147, pp. 1059--1064.

\bibitem{JKB94}
N. L. Johnson, S. Kotz, N. Balakrishnan, \emph{Continuous Univariate Distributions. Vol. 1}, second ed., Wiley Series in Probability and Mathematical Statistics: Applied Probability and Statistics, John Wiley \& Sons, Inc., New York, 1994.

\bibitem{Kolbig1972}
K. S. K\"olbig, On the zeros of the incomplete gamma function, \emph{Math. Comp.} \textbf{26} (1972), no. 119, pp. 751--755.

\bibitem{Lukas2014}
M. A. Lukas, Performance criteria and discrimination of extreme undersmoothing in nonparametric regression, \emph{J. Stat. Plan. Inference} \textbf{153} (2014), pp. 56--74.

\bibitem{Mahler1930}
K. Mahler, \"Uber die Nullstellen der unvollst\"andigen Gammafunktionen, \emph{Rend. del Circ. Matem. Palermo} \textbf{54} (1930), pp. 1--41.

\bibitem{MS16}
C.~Mitschi, D.~Sauzin, \emph{Divergent series, summability and resurgence. {I}}, Lect. Notes Math. \textbf{2153},
Springer, 2016.

\bibitem{Navarra2010}
A. Navarra, C. M. Pinotti, V. Ravelomanana, F. Betti Sorbelli, R. Ciotti, Cooperative training for high density sensor and actor networks, \emph{J. Sel. Areas Commun.} \textbf{28} (2010), no. 5, pp. 753--763.

\bibitem{Nemes13}
G. Nemes, An explicit formula for the coefficients in Laplace's method, \emph{Constr. Approx.} \textbf{38} (2013), no. 3, pp. 471--487.

\bibitem{Nemes15b}
G. Nemes, The resurgence properties of the incomplete gamma function II, \emph{Stud. Appl. Math.} \textbf{135} (2015), no. 1, pp. 86--116.

\bibitem{Nemes16}
G. Nemes, The resurgence properties of the incomplete gamma function, I, \emph{Anal. Appl. (Singap.)} \textbf{14} (2016), no. 5, pp. 631--677.

\bibitem{NIST:DLMF}
\emph{NIST Digital Library of Mathematical Functions}. \href{http://dlmf.nist.gov/}{http://dlmf.nist.gov/}, Release 1.0.19 of 2018-06-22. F. W. J. Olver, A. B. Olde Daalhuis, D. W. Lozier, B. I. Schneider, R. F. Boisvert, C. W. Clark, B. R. Miller, and B. V. Saunders, eds.

\bibitem{OD98a}
A. B. Olde Daalhuis, On the resurgence properties of the uniform asymptotic expansion of the incomplete gamma function, \emph{Methods Appl. Anal.} \textbf{5} (1998), no. 4, pp. 425--438.

\bibitem{Pagurova1956}
V. I. Pagurova, An asymptotic formula for the incomplete gamma function, \emph{U.S.S.R. Comput. Math. and Math. Phys.} \textbf{5} (1965), no. 1, pp. 162--166.

\bibitem{Paillard2008}
G. Paillard, V. Ravelomanana, Limit theorems for degree of coverage and lifetime in large sensor networks, \emph{IEEE INFOCOM 2008 - The 27th Conference on Computer Communications}, Phoenix, Arizona, USA, 2008, pp. 106--110.

\bibitem{Palffy2017}
P. Palffy-Muhoray, E. G. Virga, X. Zheng, Onsager's missing steps retraced, \emph{J. Phys.: Condens. Matter} \textbf{29} (2017), no. 47, Article 475102, 13 pp.

\bibitem{Paris2002a}
R. B. Paris, Error bounds for the uniform asymptotic expansion of the incomplete gamma function, \emph{J. Comput. Appl. Math.} \textbf{147} (2002), no. 1, pp. 215--231.

\bibitem{Paris2002b}
R. B. Paris, A uniform asymptotic expansion for the incomplete gamma function, \emph{J. Comput. Appl. Math.} \textbf{148} (2002), no. 2, pp. 323--339.

\bibitem{Paris2018}
R. B. Paris, A uniform asymptotic expansion for the incomplete gamma functions revisited, preprint, \href{https://arxiv.org/abs/1611.00548}{arXiv:1611.00548}

\bibitem{Ravelomanana2004}
V. Ravelomanana, Extremal properties of three-dimensional sensor networks with applications, \emph{IEEE Trans. Mobile Comput.} \textbf{3} (2004), no. 3, pp. 246--257.

\bibitem{Temme1979}
N. M. Temme, The asymptotic expansion of the incomplete gamma functions, \emph{SIAM J. Math. Anal.} \textbf{10} (1979), no. 4, pp. 757--766.

\bibitem{Temme1982}
N. M. Temme, The uniform asymptotic expansion of a class of integrals related to cumulative distribution functions, \emph{SIAM J. Math. Anal.} \textbf{13} (1982), no. 2, pp. 239--253.

\bibitem{Temme1992}
N. M. Temme, Asymptotic inversion of incomplete gamma functions, \emph{Math. Comp.} \textbf{58} (1992), no. 198, pp. 755--764.

\bibitem{Temme1994}
N. M. Temme, Computational aspects of incomplete gamma functions with large complex parameters, \emph{Approximation and Computation: A Festschrift in Honor of Walter Gautschi} (R. V. M. Zahar, ed.), Internat. Ser. of Numer. Math., Vol. 119, 1994, pp. 551--562.

\bibitem{Temme1996b}
N. M. Temme, Uniform asymptotics for the incomplete gamma functions starting from negative values of the parameters, \emph{Methods Appl. Anal.} \textbf{3} (1996), no. 3, pp. 335--344.

\bibitem{Thompson2012}
I. Thompson, A note on the real zeros of the incomplete gamma function, \emph{Integral Transforms Spec. Funct.} \textbf{23} (2012), no. 6, pp. 445--453.

\bibitem{Trailovic2004}
L. Trailovi\'c,  L.Y. Pao, Computing budget allocation for efficient ranking and selection of variances with application to target tracking algorithms, \emph{IEEE Trans. Automat. Contr.} \textbf{49} (2004), no. 1, pp. 59--67.

\bibitem{TW14}
L.~N. Trefethen, J.~A.~C. Weideman, The exponentially convergent trapezoidal rule, \emph{SIAM Rev.} \textbf{56} (2014), pp. 385--458.

\bibitem{Tricomi1950}
F. G. Tricomi, Asymptotische Eigenschaften der unvollst\"andigen Gammafunktion, \emph{Math. Z.} \textbf{53} (1950), pp. 136--148.

\end{thebibliography}
\end{document}